\documentclass[11pt]{article}
\usepackage{amsmath,amsfonts,latexsym,amssymb,amscd, graphicx,pictexwd, mathrsfs, dsfont, color,mathtools}
\usepackage[title]{appendix}

\newcommand{\F}{{\mathbf F}}

\newcommand{\LL}{{\mathbb L}}
\newcommand{\CC}{{\mathbf C}}

\renewcommand{\SS}{{\mathbf S}}
\newcommand{\K}{{\mathbf K}}
\newcommand{\I}{{\mathbf I}}

\renewcommand{\P}{{\mathbb P}}
\newcommand{\E}{{\mathbb E}}
\newcommand{\bF}{{\mathbb F}}
\newcommand{\bM}{{\mathbb M}}

\newcommand{\R}{{\mathbb R}}

\newcommand{\D}{{\mathbb D}}
\newcommand{\Q}{{\mathbb Q}}

\newcommand{\ZZ}{{\mathbb Z}}

\newcommand{\cA}{{\mathcal A}}
\newcommand{\cL}{{\mathcal L}}
\newcommand{\cF}{{\mathcal F}}

\newcommand{\cN}{{\mathcal N}}

\newcommand{\fh}{{\mathfrak h}}
\newcommand{\fb}{{\mathfrak b}}

\newcommand{\ff}{{\mathfrak f}}
\newcommand{\fd}{{\mathfrak d}}

\newcommand{\sfUC}{{\mathsf{UC}}}

\newcommand{\supp}{{\mathsf{supp}}}

\newcommand{\1}{{\mathds{1}}}

\newcommand{\argmax}{\mathop{\rm arg\,max}}

\newcommand{\Wass}{{\rm Wass}}
\newcommand{\Lip}{{\rm Lip}}

\newcommand{\Var}{{\rm Var\,}}
\newcommand{\Cov}{{\rm Cov}}

\newtheorem{thm}{Theorem}

\newtheorem{lem}{Lemma}

\newtheorem{rem}{Remark}

\numberwithin{equation}{section}

\begin{document}
\title{Integration by Parts and the KPZ Two-Point Function}
\author{Leandro P. R. Pimentel}
\date{\today}
\maketitle

\begin{abstract}  
In this article we consider the KPZ fixed point starting from a two-sided Brownian motion with an arbitrary diffusion coefficient. We apply the integration by parts formula from Malliavin calculus to establish a key relation between the two-point (correlation) function of the spatial derivative process and the location of the maximum of an Airy process plus Brownian motion minus a parabola. Integration by parts also allows us to deduce the density of this location in terms of the second derivative of the variance of the KPZ fixed point. In the stationary regime, we find the same density related to limit fluctuations of a second-class particle. We further develop an adaptation of Stein's method that implies asymptotic independence of the spatial derivative process from the initial data. 
\end{abstract}

\section{Introduction}
Growth models in the one dimensional Kardar-Parisi-Zhang (KPZ) universality class are usually described by a growing  interface represented by height function $h(x,t)$ at time $t\geq 0$, over a one-dimensional substrate $x\in\R$, whose evolution undergoes a stochastic local dynamics subject to three key features: smoothing, slope dependent growth speed and space-time locally correlated noise. The canonical example is the KPZ equation \cite{KPZ} $\partial_t h=\partial^2_xh+(\partial_xh)^2+\xi$, where $\xi$ is a space-time white noise. The scaling behaviour of this equation should be the same as that of the models that share these three features, and the long time limit distributions should be universal within certain geometry dependent subclasses that are encoded by the scaling properties of the initial growth profile \cite{Co}. Illustrations of natural growth phenomena within this universality class include turbulent liquid crystals, bacteria colony growth and paper wetting \cite{IwFuTa,Ta}. All models in the KPZ universality class, under the $1:2:3$ scaling transformation  $h(x,t)\mapsto c_1\epsilon^1 h(c_2\epsilon^{-2} x,\epsilon^{-3} t)-C_\epsilon t$, are conjectured to converge   to a universal fluctuating field $\fh_t(x)$, as $\epsilon\searrow 0$, called the KPZ fixed point \cite{CoQuRe}. The prefactor $C_\epsilon$ is the macroscopic speed at $x=0$ and the constants $c_1$ and $c_2$ may depend on the distributional details of the model. The geometry dependent subclasses are then determined by the initial profile $\fh_0(x)$.  
\medskip

In the last two decades a great progress was made by considering stochastic integrable growth models, notably the totally asymmetric simple exclusion process (TASEP) \cite{BaFePe0,BoFePrSa,Jo1,MaQuRe,Sa} and the polynuclear  growth model (PNG) \cite{BaDeJo,BaRa,Jo2,PrSp1}, that lead to a detailed description of the Markov interface evolution $(\fh_t(x)\,,\,x\in\R)_{t\geq 0}$ through the calculation of its functional transition probabilities as a Fredholm determinant formula. For the Cole-Hopf solution of the KPZ equation the one-point marginal limit fluctuations were proved in \cite{AmCoQu,BoCoFeVe}. The state space for the KPZ fixed point is the collection $\sfUC$ of upper semicontinuous generalised functions satisfying a linear growth control. Although the initial data in $\sfUC$ can be very irregular, for every $t>0$ the process takes values in real valued functions that look locally like Brownian motion \cite{CoHa,MaQuRe,Pi1,SoVi} and, for fixed $x\in \R$, the time evolution is locally $1/3-$ H\"older continuous \cite{MaQuRe}. Three renowned examples are as follows: The Airy$_2$ process (narrow wedge initial profile),
$$\fh_1(x)+x^2=\cA_2(x)\,,\mbox{ where $\fh_0(x)=-\infty$ for all $x\neq 0$ and $\fh_0(0)=0$}\,;$$ 
The Airy$_1$ process  (flat initial profile),
$$\fh_1(x)=2^{1/3}\cA_1(2^{-2/3}x)\,,\mbox{ where $\fh_0(x)=0$ for all $x\in\R$}\,;$$
The Airy$_{\rm stat}$ process (Brownian initial profile),
$$\fh_1(x)=\cA_{\rm stat}(x)\,,\mbox{ where }\fh_0\equiv\mbox{ two-sided Brownian motion with $\Var(\fh_0(1))=2$}\,.$$
The two-sided Brownian motion with $\Var(\fh_0(1))=2$ plays an important role since it gives the stationary KPZ fixed point in the sense that the distribution of $(\fh_t(x)-\fh_t(0)\,,\,x\in\R)$ remains the same at all times. These three subclasses are characterized by the following marginal distributions at $x=0$: the GUE Tracy-Widom \cite{TrWi1}, the GOE Tracy-Widom \cite{TrWi2} and and Baik-Rains distributions \cite{BaRa}, respectively. For general initial data, where $\fh_0=\fh$ is fixed, the one-point marginal 
$$\F_{\fh}(x,t,r):=\P_\fh\left[\fh_t(x)\leq r\right]\,,$$
are differentiable with respect to $(x,t,r)$ and $\partial^2_r\log \F_\fh$ solves the KP-II equation \cite{QuRe}.
\medskip

The KPZ fixed point also has a variational description \cite{CoLiWa,CoQuRe,MaQuRe}
\begin{equation}\label{KPZMax}
\fh_1(x)\stackrel{dist.}{=}\max_{z\in\R}\left\{\fh(z)+\cA_2(z)-(x-z)^2\right\}\,,
\end{equation}
that points the importance of the Airy$_2$ process. The variational formula \eqref{KPZMax} is valid in a broader sense and the KPZ fixed point has an alternative description as a stochastic flow constructed from the directed landscape $\cL(z,s;x,t)$ \cite{DaOrVi}, a random continuous four dimensional field that can be seen as a metric between the space time points $(z,s)$ and $(x,t)$, with $s<t$. In this article we apply the integration by parts formula from Malliavin calculus to analyse the two-point (correlation) function of the distribution valued spatial derivative process $\partial_x\fh_t$, where at time zero we have white noise with strength $\beta^2$. This process can be thought of as the stochastic Burgers fixed point, presumably the scaling limit of the stochastic Burgers equation $\partial_t u=\partial_x u^2+\partial^2_x u+\partial_x\xi$, that relates to the KPZ equation by the transformation $u=\partial_x h$. We show that the directed landscape provides a geometrical description of the Malliavin derivative of an observable of $\partial_x\fh_t$ in terms of the a.s. unique location of the maximum in \eqref{KPZMax} at $x=0$ \cite{Pi1}, defined as
\begin{equation}\label{ArgEq}
Z:=\argmax_{z\in\R}\left\{\fh(z)+\cA_2(z)-z^2\right\}\,.
\end{equation}
This description is combined with integration by parts to derive the two-point function, and the density of $Z$, in terms of the second derivative of the variance of $\fh_1(x)$. We go beyond the analysis of the correlation and combine the previous results with Malliavin-Stein method to prove that the joint law of observables of the system with respect to $\partial_x\fh_0$ and $\partial_x\fh_t$ is  close, in the Wasserstein metric, to the product measure induced by its marginals (asymptotic independence). Next we explain  the main results of this article with more details.      
\medskip

The analysis is restricted to an initial profile $\fh_0(x)=\beta \fb(x)$, where $\beta>0$ and the stochastic process $\fb\equiv(\fb(x)\,,\,x\in\R)$ is a standard two-sided Brownian motion. In order to not overload notation, we keep the dependency on the parameter $\beta$ implicit, and the underlying probability measure $\P=\P_\beta$ is the product measure induced by the initial profile and the transition probabilities of the KPZ fixed point. This class of distributions labeled by $\beta\geq 0$ was also considered in \cite{ChFeSp}, where you can see the plot of the probability densities of $\fh_1(0)$ obtained by TASEP Monte Carlo simulations (with $\beta^2=2\sigma^2$ there). Recall that the time stationary regime is given by $\Var\fh_0(1)=\beta^2=2$. The observable of the system at time $t\geq 0$ is defined as
$$X^{\phi}_t:=\int_\R\phi \partial_x\fh_tdx\,,$$
where $\phi:\R\to\R$ is a given deterministic test function. Later in the text we provide a precise meaning of the integral, and what kind of test functions we are considering. Since $\fh_0$ is assumed to be a Brownian motion, $X^{\phi}_0$ is defined as the Wiener integral of $\phi$ with respect to $\partial_x\fh_0dx=\beta d\fb $, but for $t>0$ we need to be more careful in the definition of $X^{\phi}_t$. Observe that $X_0^\phi$ is Gaussian random variable with mean zero and variance $\beta^2  \|\phi \|^2_{\LL^2(\R)}$, where $\|\phi \|^2_{\LL^2(\R)}$ denotes the $\LL^2$-norm of $\phi\in\LL^2(\R)$. 
\medskip

The KPZ two-point function is given by the cross correlation of the differentials $\partial_x\fh_0(u)$ and $\partial_x\fh_t(v)$, which is expected to be a nonnegative function\footnote{This claim is motivated by the study of the totally asymmetric simple exclusion process, where the two-point function is proportional to the probability function of a second-class particle \cite{PrSp2}.}  $C_\beta$ of $z=v-u$ (by space stationarity of the initial data) and $t$ such that 
\begin{equation}\label{CrossCorr0}
\Cov(X_0^{\phi_1},X_t^{\phi_2})=\E\left[X_0^{\phi_1}X_t^{\phi_2}\right]=\int_{\R^2}\phi_1(u)\phi_2(v)C_\beta(v-u,t)dudv\,
\end{equation}
Notice that, if \eqref{CrossCorr0} is true then by a simple change of variable and Fubini's theorem, it can be rewriten as 
$$\E\left[X_0^{\phi_1}X_t^{\phi_2}\right]=\int_{\R^2}\phi_1(u)\phi_2(v)C_\beta(v-u,t)dudv=\int_\R \phi_1\star \phi_2(z)C_\beta(z,t)dz\,,$$
in terms of the cross correlation of the test functions $\phi_1$ and $\phi_2$, defined as  
\begin{equation}\label{CrossCorr}
\phi_1\star\phi_2(z):=\int_\R\phi_1(u)\phi_2(u+z)du\,.
\end{equation}
The first result in order to prove \eqref{CrossCorr0} and to determine $C_\beta$ is the following (Theorem \ref{CovArg}): 
\begin{equation}\label{CovArgEq}
\E\left[X_0^{\phi_1} X_t^{\phi_2}\right]=\beta^2\E\left[\phi_1\star\phi_2\left(t^{2/3}Z\right)\right]\,,
\end{equation}
where $Z$ is given by \eqref{ArgEq}. The proof of \eqref{CovArgEq} uses the integration by parts formula from Malliavin calculus, which naturally appears to express the left hand side of \eqref{CovArgEq} as the expected value of the $\LL^2(\R)$ inner product between the Malliavin derivative of $X^{\phi_2}_t$ with respect to the initial data $\fb$ and the test function ${\phi_1}$. 
\medskip

Formula \eqref{CovArgEq} indicates the relation between the KPZ two-point function and the density of $Z$, and we can actually use integration by parts again to compute the distribution of $Z$ in terms of the variance of $\fh_1(x)$ as follows (Theorem \ref{VarDen}). The KPZ scaling function $g_\beta$ is defined as   
$$g_\beta(x):=\Var\fh_1(x)\,,$$
where the variance is computed with respect to the product measure $\P=\P_\beta$ induced by the random initial profile $\beta\fb$ and the transition probabilities of the KPZ fixed point. We denote the distribution function of $Z$ by 
$$\bF_\beta(x):=\P\left[Z\leq x\right]=\P_\beta\left[Z\leq x\right]\,.$$ 
We will prove that $g_\beta$ is differentiable with respect to $x$ and that (Theorem \ref{VarDen})
\begin{equation}\label{Eq1VarDen}
g'_\beta(x)=\beta^2\left(2\bF_\beta(x)-1\right)\,.
\end{equation}
Differentiability of $g_\beta$ also follows from differentiability of the kernel in the Fredholm determinant formula for the probability law of $\fh_1(x)$ \cite{MaQuRe,QuRe} (plus some uniform estimates to differentiate under the integral sign). The function $g_\beta$ is twice differentiable with respect to $x$  and, as a corollary of \eqref{CovArgEq} and \eqref{Eq1VarDen}, $Z$ has the density  
\begin{equation}\label{Eq2VarDen}
f_\beta(x)=\frac{g''_\beta(x)}{2\beta^2}\,,
\end{equation} 
which finaly shows that 
\begin{equation}\label{Eq3VarDen}
\E\left[X_0^{\phi_1} X_t^{\phi_2}\right]=\int_\R\phi_1\star\phi_2(z)\frac{g''_\beta(zt^{-2/3})}{2t^{2/3}}dz\,\,\mbox{ and }\,\,C_\beta(z,t)=\frac{g''_\beta(zt^{-2/3})}{2t^{2/3}}\,.
\end{equation}

In the stationary regime $\beta^2=2$ the KPZ scaling function is commonly denoted $g_{sc}\equiv g_{\sqrt{2}}$, and \eqref{Eq3VarDen} was already obtained by using TASEP approximations to the KPZ fixed point \cite{BaFePe,FeSp,MaQuRe,PrSp2}. The distribution of $\fh_1(x)$ has an explicit formulation in terms of the Painlev\'e II equation \cite{BaRa,PrSp2}, which is also related to the KP-II equation \cite{QuRe}, or in terms of the Tracy-Widom GUE distribution and the Airy kernel \cite{BaFePe,FeSp}. In \cite{PrSp2} one can see the numerical method developed to compute $f_{KPZ}\equiv f_{\sqrt{2}}= g''_{sc}/4$, as the plot of its graph. As it was point out in \cite{FeSp,PrSp2}, the same density \eqref{Eq2VarDen} appears as the limit fluctuations of a second class particle in the PNG and TASEP stationary regimes, which is consistent with the well known duality between maximizers and second-class particles \cite{BaCaSe,CaGr}. In the physics literature, the relation between the density of the location of the maximum in the stationary regime and $g''_{sc}/4$ was predicted using Bethe ansatz calculations \cite{Le}, while in \cite{MaTh} it was based on a linear response method applied to the stationary stochastic Burgers equation \cite{ImSa}, which resembles in some aspects the computation of the Malliavin derivative at a fixed direction (compare equation (8) in \cite{MaTh} with \eqref{IntPart} in next section). KPZ correlations starting from Brownian profiles were also tested in experimental systems of growing liquid crystal turbulence \cite{IwFuTa}.   
\medskip

Some of the novelties in this article are: the rigorous deduction of \eqref{Eq3VarDen} in terms of the variance of $\fh_1(x)$ for every value of $\beta>0$, which is related to the predictions raised in  \cite{ChFeSp} about the behaviour of the two-point function; the explicit relation \eqref{Eq2VarDen} with the density of $Z$; the method of proof by means of Malliavin calculus applied to the directed landscape formulation of the KPZ fixed point. For flat profile $\fh_0\equiv 0$, which corresponds to $\beta=0$, the distribution of $Z$ was computed in \cite{MoQuRe} by a direct calculation using a Fredholm determinant formula for the probability that $\cA_2(z)\leq a(z)$ on a finite interval \cite{CoQuRe0}, where $a(z)$ was essentially a parabola. In the same $\beta=0$ regime, a different formula for the density of $Z$ was obtained in \cite{BaLiSc,Sc}. In the Brownian case the same type of determinant calculation faces the problem that $a(z)$ is going to be a parabola plus a sample of a Brownian motion, and then one has to integrate out the resulting formula with respect to Brownian motion, which seems to lead us to a problem with no way out. We note that by space stationarity of the Airy$_1$ process $g_0(x)=g_0(0)$ for all $x\in\R$. On the other hand, since $\bF_\beta\to\bF_0$, as $\beta\searrow 0$, by \eqref{Eq1VarDen}, one has that  
\begin{equation}\label{flat}
\frac{\partial_\beta^2 g'_0(x)}{2}=\lim_{\beta\searrow 0}\beta^{-2}g_\beta'(x)=2\bF_0(x)-1\,.
\end{equation}
This points out an alternative way to compute the distribution of $Z$ in the flat case, as soon as one can obtain an explicit  formula for the left hand side of \eqref{flat}, and a natural candidate for $\bF_0$ arises by taking the second derivative of $g'_\beta(x)$ with respect to $\beta$ and eveluating at $\beta=0$. Another interesting aspect is the $\beta\nearrow\infty$ regime \cite{ChFeSp}, where the Airy$_2$ process becomes irrelevant and 
$$\lim_{\beta\nearrow\infty}\beta^{-2/3}Z\stackrel{dist.}{=}\argmax_{u\in\R}\left\{\fb(u)-u^2\right\}\,,$$
which has the well known Chernoff's distribution $\bF_{\rm Ch}$ \cite{GrWe}. Thus, by \eqref{Eq1VarDen}, 
\begin{equation}\label{Chernoff}
\lim_{\beta\nearrow \infty}\frac{g_\beta'(\beta^{2/3}x)}{\beta^{2}}=2\bF_{\rm Ch}(x)-1\,,
\end{equation}
and one might also expect that
$$\lim_{\beta\nearrow \infty}\frac{g_\beta''(\beta^{2/3}x)}{2\beta^{4/3}}=f_{\rm Ch}(x)\,,$$
where $f_{\rm Ch}$ is the density of $\bF_{\rm Ch}$.

\medskip

After studying the cross correlation between $X_0^{\phi_1}$ and $X_t^{\phi_2}$, we turn to a more delicate issue related to quantifying the distance between the joint law $\theta_t=\P_{X_0^{\phi_1},X_t^{\phi_2}}$ of the random vector $(X_0^{\phi_1},X_t^{\phi_2})$, and the product measure $\eta_t=\P_{X_0^{\phi_1}}\otimes\P_{X_t^{\phi_2}}$ induced by its  marginals. To reach that goal we develop a simple adaptation of Malliavin-Stein method, which allows us to prove the following upper bound (Theorem \ref{AsyInd}):
\begin{eqnarray}
\nonumber\Wass\left(\eta_t,\theta_t\right)&\leq&\frac{\beta}{\|\phi_1\|_{\LL^2(\R)}}\sqrt{\frac{\pi}{2}} \E\left[\left |\psi_1\right |\star\left |\phi'_2\right |\left(t^{2/3}Z\right)\right]\\
\label{EqAsyInd}&=&\frac{1}{\beta\|\phi_1\|_{\LL^2(\R)}}\sqrt{\frac{\pi}{2}} \int_\R\left |\psi_1\right|\star\left |\phi'_2\right|(z)\frac{g''_\beta(zt^{-2/3})}{2t^{2/3}} dz\,,
\end{eqnarray}
where $\Wass(\eta_t ,\theta_t)$ denotes the Wasserstein distance between the probability measures $\eta_t$ and $\theta_t$, and $\psi_1'=\phi_1$ with $\psi_1(0)=0$.

\medskip
By \eqref{EqAsyInd}, the the distance to independence exhibited by the joint law of $X_0^{\phi_1}$ and $X_t^{\phi_2}$ scales as $t^{2/3}$, and the limit behaviour is connected to 
\begin{equation}\label{BoundDen}
\lim_{t\to\infty}\int_\R|\psi_1|\star|\phi'_2|(z)\frac{g''_\beta(zt^{-2/3})}{2} dz\,=\,\frac{g_\beta''(0)}{2}\int_\R |\psi_1|(u)du \int_\R |\phi'_2|(v)dv\,.
\end{equation}
However to justify \eqref{BoundDen}, one needs more information on $g''_\beta$. This function should be bounded by $g_\beta''(0)$ for all $\beta >0$, which is sufficient to obtain \eqref{BoundDen} (by dominated convergence). For $\beta=\sqrt{2}$ this is known \cite{PrSp2}, and a numerical computation shows that $g''_{sc}(0)\approx 2.16$. For $\beta\searrow 0$ one can use that $g''_\beta(0)\approx 2\beta^2 f_0(0)$ \cite{MoQuRe} while for $\beta\nearrow\infty$, $g''_{\beta}(0)\approx 2\beta^{4/3}f_{\rm Ch}(0)$ \cite{GrWe}. 
\medskip

The relation between Malliavin calculus and Stein's method is a current research topic in the field of stochastic partial differential equations, where it can be used to prove a central limit theorem for the spatial integral of a solution \cite{HuNuViZh,Nu}. To the best of our knowledge, it is the first time that this two subjects are combined to obtain asymptotic independence from the initial data. The ideas presented in Stein's methods are adequately general to be able to apply to approximations by distributions other than the normal, such as Poisson, binomial or exponential \cite{Ro}, and it can be used to prove asymptotic independence for other Markov processes as well, such as queues in tandem, particle systems and solutions of stochastic differential equations, which is left for future works.

\noindent\paragraph{\bf Organization} 
In Section \ref{KPZ} we give a more detailed introduction of the KPZ fixed point and state the main results: Theorem \ref{CovArg}, Theorem \ref{VarDen} and Theorem \ref{AsyInd}. In Section \ref{L2} we prove some $\LL^2$ estimates for the KPZ fixed point, and in Section \ref{MalCal} we introduce the basic tools from Malliavin calculus and  prove Theorem \ref{CovArg} and Theorem \ref{VarDen}. In Section \ref{SteMet} we develop an adaptation of Stein's method, having in mind the product measure as our target measure, and prove Theorem \ref{AsyInd}.
\medskip

\noindent\paragraph{\bf Acknowledgement} The author would like to thank Patrik Ferrari, Jeremy Quastel for useful comments and enlightening discussions concerning this subject, and to thank Daniel Remenik for point me out the differentiability of $g_\beta$ from the Fredholm determinant formula for the KPZ fixed point \cite{QuRe}. Much of this work was developed during the XXIII Brazilian School of Probability, and highly inspired by the Malliavin's Calculus classes given by D. Nualart \cite{Nu}, for which the author is very grateful. This research was supported in part  by the National Council of Scientific Researches (CNPQ, Brazil) grant 305356/2019-4.

\section{The KPZ Fixed Point}\label{KPZ}
 
Let $\sfUC$ denote the space of functions $\ff:\R\to\R\cup\{ -\infty\}$ such that: (i) $\limsup_{x\to y}\ff(x)\leq \ff(y)$ (upper semicontinuity); (ii) $\ff(x)\leq C_1|x|+C_2$ for all $x\in\R$, for some $C_1,C_2<\infty$; (iii) $\ff(x)>-\infty$ for some $x\in\R$. The state space $\sfUC$ can be endowed with the topology of local convergence turning it into a Polish space (Section 3.1 \cite{MaQuRe}), such that the collection composed by cylindrical subsets of $\sfUC$,    
$$\mathrm{Cy}(\vec{x},\vec{a}):=\Big\{\ff\in\sfUC\,:\,\ff(x_1)\leq a_1,\dots,\ff(x_m)\leq a_m\Big\}\,\mbox{ for }\vec{x},\vec{a}\in\R^m\,,$$
is a generating sub-algebra for the Borel $\sigma$-algebra over $\sfUC$. The KPZ fixed point $\left(\fh_t(\cdot)\,,\,t\geq 0\right)$, with $\fh_0=\fh\in\sfUC$, is the unique time homogeneous  Markov process taking values in $\sfUC$ with transition probabilities given by
\begin{equation}\label{DefKPZ}
\P_\fh\Big[\fh_t\in \mathrm{Cy}(\vec{x},\vec{a})\Big]=\det\left[\I-\K^{\fh}_{t,\vec{x},\vec{a}}\right]_{\LL^2(\{x_1,\dots,x_m\}\times\R)}\,,
\end{equation}
when restricted to the sub-algebra composed by cylindrical subsets. This process was introduced by Matetski, Quastel and Remenik (Definition 3.12 in \cite{MaQuRe}) to describe the limit fluctuations of the rescaled height function associated to the TASEP, started from an initial data for which the diffusive scaling limit is given by $\fh$. On the right hand side of \eqref{DefKPZ} we have a Fredholm determinant of the integral  operator $\K^{\fh}_{t,\vec{x},\vec{a}}$, whose definition we address to \cite{MaQuRe} ($\I$ is the identity operator), where we have the counting measure on $\{x_1,\dots,x_m\}$ and the Lebesgue measure on $\R$. As mentioned in the introduction, from this formula one can recover several of the classical Airy processes by starting with special profiles for which the respective operators are explicit (see Section 4.4 of \cite{MaQuRe}). One of the central features of the KPZ fixed point is the so called 1:2:3 scaling invariance:
\begin{equation}\label{123}
S_{\gamma^{-1}}\fh_{\gamma^{-3}t}(\cdot;S_\gamma\fh)\stackrel{dist.}{=}\fh_{t}(\cdot;\fh)\,,\,\mbox{ where }S_\gamma\ff(x):=\gamma^{-1}\ff(\gamma^2 x)\,.
\end{equation}

This Markov process has an alternative description in terms of the a variational formula initially introduced by Corwin, Quastel and Remenik \cite{CoQuRe}, and then rigorously constructed by Dauvergne, Ortmann and Vir\'ag \cite{DaOrVi} in terms of the directed landscape, the unique four-dimensional continuous random field $\cL:\R_\uparrow^4\to\R$, where $\R_\uparrow^4:=\left\{(z,s;x,t)\,:\,s<t\,\mbox{ and }\,z,x\in\R\right\}$, that satisfies the following properties\footnote{In \cite{DaVi,NiQuRe} it was proved that the variational formula \eqref{EvoDef} describes the KPZ fixed fixed point defined through the transitions \eqref{DefKPZ}, and also convergence of several integrable models to the KPZ fixed point.}.
\begin{itemize}
\item Independent increments: if $\{(t_i,t_i+s_i)\,:\,i=1,\dots,k\}$ is a  collection of disjont intervals then $\{\cL(\cdot,t_i;\cdot,t_i+s_i)\,:\,i=1,\dots,k\}$ is a collection of independent two-dimensional random fields. 
\item Metric composition: almost surely
\begin{equation*}\label{MetComp}
\cL(x,r;y,t)=\max_{z\in\R}\Big\{\cL(x,r;z,s)+\cL(z,s;y,t)\Big\}\,,\,\forall\,(x,r;y,t)\in\R_{\uparrow}^4\mbox{ and }s\in(r,t)\,.
\end{equation*}
\item Airy sheets marginals: for fixed time $t\in\R$ and $s>0$
\begin{equation*}
\left\{\cL(x,t;y,t+s^3)\,:\,(x,y)\in\R^2\right\}\stackrel{dist.}{=}\left\{s\cL(x/s^2,y/s^2)\,:\,(x,y)\in\R^2\right\}\,,
\end{equation*}
where $\cL(x,y):=\cA(x,y)-(x-y)^2$ and $\cA:\R^2\to\R$ is a random stationary and symmetric field, namely the Airy sheet, that is uniquely determined as a functional of the Airy line ensemble \cite{CoHa,DaOrVi}. Furthermore, for fixed $y\in\R$, 
$$\left\{\cA(x,y)\,:\,x\in\R\right\}\stackrel{dist.}{=}\left\{\cA_2(x)\,:\,x\in\R\right\}\,.$$ 
\end{itemize}
Due to the parabolic drift towards $-\infty$, a.s. for all $s<t$ and $x\in\R$, the random function $z\in\R\mapsto\fh(z)+\cL(z,s;x,t)$ attains its maximum on a compact set and, due to metric composition, the process 
\begin{equation}\label{EvoDef}
\fh_{s,t}(x;\fh):=\max_{z\in\R}\Big\{\fh(z)+\cL(z,s;x,t)\Big\}\,,
\end{equation} 
defines a time homogeneous  Markov evolution acting on $\sfUC$, with $\fh_{0,0}=\fh$ and transition probabilities determined by \eqref{DefKPZ} \cite{NiQuRe}. By \eqref{EvoDef}, 
$$\cA(x,y) = \fh_{0,1}(y;\fd_x)+(x-y)^2\,,\,\mbox{ where }\,\fd_x(z)=\left\{\begin{array}{ll}0 & \mbox{ for } z=x\\-\infty &\mbox{ for } z\neq x\,.\end{array}\right.$$
The set  
$$\argmax_{z\in\R}\Big\{\fh(z) +\cL(z,s;x,t)\Big\}:=\Big\{z\in\R\,:\,\fh(z) +\cL(z,s;x,t)=\fh_{s,t}(x;\fh)\Big\}\,,$$
is compact and we also consider the process defined by the rightmost location of the maximum  
\begin{equation}\label{ArgMax}
Z_{s,t}(x;\fh):=\max\argmax_{z\in\R}\Big\{\fh(z)+\cL(z,s;x,t)\Big\}\,.
\end{equation}
For $s=0$, fixed $t>0$ and $x\in\R$, a.s. the maximum is attained at a unique location (see Section 3.1 and Proposition 5 in \cite{Pi1}), however this is not true simultaneously for all $x\in\R$ \cite{CoHaHeMa}. 
\medskip

In what follows, given a measure space $\bM$, we denote $\|\cdot\|_{\LL^p(\bM)}$ the usual $\LL^p(\bM)$ norm and $\langle\cdot,\cdot\rangle_{\LL^2(\bM)}$ the usual $\LL^2(\bM)$ inner product. When $\bM=\R$ we are always considering the Lebesgue measure over $\R$ endowed with the Borel $\sigma$-algebra.  We work with the following additional function spaces: the space $\CC_b^1(\R)$ of all continuously differentiable real valued functions on $\R$ of bounded support; the space $\SS(\R)$ of step real valued functions on $\R$ of bounded support, that is there exist $n\geq 1$, $c_i\in\R$ and $x_{i-1}\leq x_i$ for $i=1,\dots,n$ such that $\phi(x)=c_i$ for $x\in(x_{i-1},x_i]$ and $\phi(x)=0$ if $x\leq x_0$ or $x>x_n$. Let $\fb$ denote a standard two-sided Brownian motion $(\fb(x)\,,\,x\in\R)$, that can be obtained by taking two independent standard  Brownian motions $(\fb^+(x)\,,\,x\geq 0)$ and $(\fb^-(x)\,,\,x\geq 0)$ starting at $0$, and defining $\fb(x):=\fb^+(x)$ for $x\geq 0$ and $\fb(x):=\fb^-(-x)$ for $x< 0$. We work with an initial profile $\fh_0\equiv \beta\fb$ where $\beta>0$ is a fixed parameter. By assumption, the Brownian motion $\fb$ and the directed landscape $\cL$ are independent processes.   
\medskip

As mentioned before, another important symmetry of the KPZ fixed point is related to time stationarity (up to a vertical shift): if $\fh_0=\sqrt{2}\fb$ then  
\begin{equation}\label{stat_1}
\Delta \fh_t(\cdot)\stackrel{dist.}{=}\fh_0(\cdot)\,,\,\mbox{ for all }\,t\geq 0\,,
\end{equation}
where $\Delta\ff(x):=\ff(x)-\ff(0)$ for $x\in\R$. To keep notation as simple as possible, from now on we denote 
$$\fh_t(\cdot)\equiv\fh_{0,t}(\cdot;\beta\fb)\mbox{ and }Z_t(\cdot)\equiv Z_{0,t}(\cdot;\beta\fb)\,,$$
and keep the dependence on $\beta$ implicit. Note that $Z_1(0)\stackrel{dist.}{=}Z$ as in \eqref{ArgEq}. For $t>0$ and $\beta\neq \sqrt{2}$ we consider the following integrals with respect to $\partial_x\fh_t$: 
$$\phi\in\CC^1_b(\R)\,\mapsto\,\int_\R\phi \partial_x\fh_tdx:=-\int_\R\phi'(x)\fh_t(x)dx\,,$$
and
$$\phi\in\SS(\R)\,\mapsto\,\int_\R\phi\partial_x\fh_tdx:=\sum_{i=1}^n c_i\left(\fh_t(x_i)-\fh_t(x_{i-1})\right)\,$$
where $\phi(x)=c_i$ for $x\in(x_{i-1},x_i]$ and $\phi(x)=0$ if $x\leq x_0$ or $x>x_n$. For $\beta=\sqrt{2}$, we use time stationarity \eqref{stat_1} and take the Wiener integral of $\phi\in\LL^2(\R)$ with respect to $\partial_x\fh_t$. Recall that we are considering the following observables,
\begin{equation}\label{DefObs}
X^{\phi_1}_0:=\int_\R\phi_1 \partial_x\fh_0dx\,\mbox{ and }\,X^{\phi_2}_t:=\int_\R\phi_2 \partial_x\fh_tdx\,,
\end{equation}
where $X_0^{\phi_1}$ is the Wiener integral of $\phi\in\LL^2(\R)$ with respect to $\partial_x\fh_0dx=\beta d\fb$, which is a Gaussian random variables with zero mean and variance $\beta^2\|\phi_1\|_{\LL^2(\R)}$. Recall also that the cross correlation $\phi_1\star\phi_2$ of $\phi_1$ and $\phi_2$ is defined in \eqref{CrossCorr} and the random variable $Z$ is defined in \eqref{ArgEq}. 

\begin{thm}\label{CovArg}
We have that 
\begin{equation*}
\E\left[X_0^{\phi_1} X_t^{\phi_2}\right]=\beta^2\E\left[\left(\phi_1\star\phi_2\right)\left(t^{2/3}Z\right)\right]\,.
\end{equation*}
\end{thm}

\begin{thm}\label{VarDen}
Let $g_\beta(x):=\Var\left[\fh_1(x)\right]$ and $\bF_\beta(x):=\P\left[Z\leq x\right]$. Then $g_\beta$ is differentiable and 
$$g'_\beta(x)=\beta^2\left(2\bF_\beta(x)-1\right)\,.$$
In particular, $g'_\beta$ is absolutely continuous iff $Z$ has a density $f_\beta$. In this case, we also have that 
$$f_\beta(x)=\frac{g_\beta''(x)}{2\beta^2}\,\mbox{ and }\,\E\left[X_0^{\phi_1} X_t^{\phi_2}\right]=\int_\R\left(\phi_1\star\phi_2\right)(z)\frac{g''_\beta(zt^{-2/3})}{2t^{2/3}}dz\,.$$  
\end{thm}

As it was mentioned before, twice differentiability of $g_\beta$ follows from the Fredholm determinant formula for the distribution of the KPZ fixed point \cite{MaQuRe,QuRe}. The proof of both theorems relies on the integration by parts formula from Malliavin calculus \cite{Nu} as follows. Let $\phi\in\LL^2(\R)$ and denote $W(\phi)=\int_{\R}\phi d\fb$.  If the random variable $X$ is Malliavin differentiable with respect to $\fb$ then 
$$\E\left[W(\phi)X\right]=\E\left[\langle DX,\phi\rangle_{\LL^2(\R)}\right]\,,$$
where $DX=(DX(x)\,,\,x\in\R)$ is the Malliavin derivative of $X$ with respect to $\fb$. This formula naturally leads  to 
\begin{equation}\label{IntPart}
\E\left[X_0^{\phi_1} X_t^{\phi_2}\right]=\beta\E\left[W(\phi_1) X_t^{\phi_2}\right]=\beta\E\left[\langle DX_t^{\phi_2},\phi_1\rangle_{\LL^2(\R)}\right]\,.
\end{equation}
We will show that $X^\phi_t$ is Malliavin differentiable with respect to $\fb$ and that
\begin{equation}\label{YMalDer2}
\| DX_t^\phi\|_{\LL^2(\Omega\times\R)}=\beta\|\phi\|_{\LL^2(\R)}\,\mbox{ and }\,\E\left[\langle DX_t^{\phi_2},\phi_1\rangle_{\LL^2(\R)}\right]=\beta\E\left[\left(\phi_1\star\phi_2\right)\left(t^{2/3}Z\right)\right]\,.
\end{equation}
From \eqref{IntPart} and \eqref{YMalDer2}, we  will deduce Theorem \ref{CovArg}. For a fixed realisation of $\fb$ and $\cL$, the function $x\in\R\mapsto Z_t(x)$ defined in \eqref{ArgMax} is right continuous and non-decreasing (Proposition 2 \cite{Pi1}), and by the 1:2:3 scale invariance \eqref{123}, for fixed $x\in\R$, 
$$Z_t(x)\stackrel{dist.}{=}x+t^{2/3}Z\,$$
 (Lemma \ref{ArgSc}). Relying on the coalescence property of maximal paths \cite{Pi0}, it is conjectured that the image of the function $x\in\R\mapsto Z_t(x)$ is a locally finite stationary point process. If this is true and $u_i<u_{i+1}$ denote the points of discontinuities of $Z_t$ then we also expect that  
$$DX_t^{\phi}(x)=\beta\sum_{i\in\ZZ}\phi(u_{i+1})\1_{\left(Z_t(u_i),Z_t(u_{i+1})\right]}(x)\,,$$
but we do not need this explicit form to proceed with the calculations. 
\medskip

The source of randomness comes from the initial data and the directed landscape used to construct the Markovian evolution. By assumption, these two sources are independent which allows us to use Malliavin calculus with respect to initial data $\fb$, and then integrate it out with respect to the directed landscape $\cL$. To understand from where \eqref{YMalDer2} comes from we recall that if $X$ is a functional of $\fb$, then the Malliavin derivative defines a linear (and unbounded) random operator that can be interpreted as a directional derivative \cite{Nu}: if $\phi\in\LL^2(\R)$ and $\psi(x)=\int_0^x\phi(z) dz$, with the convention that $\int_0^x\equiv -\int_x^0$ for $x<0$, then
$$\frac{d}{d\epsilon}X\left(\fb+\epsilon\psi\right)\Big |_{\epsilon=0}=\langle DX,\phi\rangle_{\LL^2(\R)}\,.$$
For $X=\fh_1(0)$ we have that (recall \eqref{KPZMax} and \eqref{ArgEq})
$$\frac{d}{d\epsilon}X\left(\fb+\epsilon\psi\right)\Big |_{\epsilon=0}=\beta\psi(Z)=\beta\int_0^Z\phi(z)dz\,\Rightarrow\, DX(z)=\left\{\begin{array}{ll}\,\,\,\,\beta\1_{(0,Z]}(z)& \mbox{ if } Z>0\\
- \beta\1_{(Z,0]}(z) &\mbox{ if } Z\leq 0\,.\end{array}\right.$$ 
In this sense, \eqref{IntPart} and \eqref{YMalDer2} can be seen as a version of the covariance formula obtained for second-class particles and exit-points in stationary TASEP and PNG, where the proofs are also based on the addition of a small perturbation to the initial profile of the system, and the computation of the rate of change as the size of the perturbation goes to zero. For instance, see the proofs of Lemma 4.6 in \cite{BaCaSe}, Theorem 2.1 in \cite{CaGr} and (2.12) in \cite{PrSp2}. To prove Theorem \ref{VarDen} we link the variance of $\fh_1(x)$ with  the covariance between $\fh_1(x)$ and $\fh_0(x)$ by a simple calculation, and compute this covariance in terms of the Malliavin derivative of $\fh_1(x)$. It is also remarkable that this simple relation between the variance and the covariance (covariance-variance reduction) was combined in \cite{DeFlOr} with tools from Malliavin calculus for concentration bounds to study aging for the stationary KPZ equation and related models.  
\medskip

After studying the covariance we turn to the problem related to quantifying the distance between the joint law $\theta_t=\P_{X_0^{\phi_1},X_t^{\phi_2}}$ and the product measure $\eta_t=\P_{X_0^{\phi_1}}\otimes\P_{X_t^{\phi_2}}$ induced by the  marginals of $\theta_t$. The Wasserstein distance between the probability measures $\eta$ and $\theta$ over $\R^2$ is defined as
\begin{equation}\label{Wass}
\Wass(\eta,\theta):=\sup\left\{\big|\,\int_{\R^2} l d\eta-\int_{\R^2} ld\theta\,\big|\,:\,l\in\Lip_1\right\}\,,
\end{equation}
where $l:\R^2\to \R$ belongs to $\Lip_C$ if 
$$\|l\|_{\Lip}:=\sup_{(x_1,x_2)\neq (y_1,y_2)}\frac{\big|l(x_1,x_2)-l(y_1,y_2)\big|}{\|(x_1,x_2)-(y_1,y_2)\|_{\R^2}}\leq C\,,$$
and $\|\cdot\|_{\R^2}$ denotes the usual euclidean norm. 

\begin{thm}\label{AsyInd}
Let $\theta_t=\P_{X_0^{\phi_1},X_t^{\phi_2}}$ denote the joint law of $(X_0^{\phi_1},X_t^{\phi_2})$ and let $\eta_t=\P_{X_0^{\phi_1}}\otimes\P_{X_t^{\phi_2}}$ denote the product measure induced by the  marginals of $\theta_t$. Then 
$$\Wass\left(\eta_t,\theta_t\right)\leq \frac{\beta}{\|\phi_1\|_{\LL^2(\R)}}\sqrt{\frac{\pi}{2}}\E\left[\left |\psi_1\right |\star\left|\phi'_2\right|\left(t^{2/3}Z\right)\right]\big |\,.$$
where $\psi_1'=\phi_1$ and $\psi_1(0)=0$.
\end{thm}

In general terms, Stein's methods is composed by two parts \cite{Ro}: (i) bound the distance between two probability measures in terms of the expectation of a certain functional of the underlying random element, that is constructed taking into account a characterizing property of the target measure (in our case, the product measure); (ii) develop techniques to bound the expectation appearing in the first part by exploring the structure provided by the random element. Let us forget about the time parameter for the moment, and denote $X\sim N(0,\sigma^2)$ if $X$ has a normal distribution with zero mean and variance $\sigma^2$. The characterizing operator of the target measure is
$$\cN f(x_1,x_2):=\sigma^2\partial_{x_1}f(x_1,x_2)-x_1f(x_1,x_2)\,,$$   
in the sense that 
$$\mbox{$X_1\sim N(0,\sigma^2)$ is independent of $X_2$ iff $\E\left[\cN f(X_1,X_2)\right]=0$}\,,$$ 
for all continuously differentiable $f:\R^2\to\R$ with bounded derivatives (it does not specify the distribution of $X_2$).  Following Stein's method, an upper bound for the Wasserstein distance between $\theta=\P_{X_1,X_2}$  and $\eta=\P_{X_1}\otimes\P_{X_2}$ is derived in terms of the unique bounded solution $f_l$ of the partial differential equation
\begin{equation}\label{SteinEq0}
\sigma^2\partial_{x_1}f(x_1,x_2)-x_1f(x_1,x_2)=l(x_1,x_2)-\E\left[l(X_1,x_2)\right]\,,
\end{equation}
where $l:\R^2\to \R$ is a $1$-Lipschitz function, in such way that the main estimate concerns the expected value of $\cN f_l(X_1,X_2)$ under the measure $\theta\equiv \P_{X_1,X_2}$ (one can actually take $l$ continuously differentiable with bounded partial derivatives). Turning back to our KPZ context, on one hand we have $X_1=X_0^{\phi_1}\sim N(0,\sigma^2)$ with $\sigma^2=\beta^2\|\phi_1\|^2_{\LL^2(\R)}$. On the other hand, the integration by parts and the chain rule for Malliavin derivatives \cite{Nu} imply that
$$\E\left[X_0^{\phi_1} f(X_0^{\phi_1},X_t^{\phi_2})\right]=\beta^2\|\phi_1\|_{\LL^2(\R)}^2\E\left[\partial_{x_1}f(X_0^{\phi_1},X_t^{\phi_2})\right]+\beta\E\left[\partial_{x_2} f(X_0^{\phi_1},X_t^{\phi_2})\langle DX_t^{\phi_2},\phi\rangle_{\LL^2(\R)}\right]\,,$$
which yields to 
\begin{equation}\label{SteinOp}
\E\left[\cN f_l(X_0^{\phi_1},X_t^{\phi_2})\right]=-\beta\E\left[\partial_{x_2} f_l(X_0^{\phi_1},X_t^{\phi_2})\langle DX_t^{\phi_2},\phi\rangle_{\LL^2(\R)}\right]\,.
\end{equation}
The analyse of \eqref{SteinEq0} shows that $\|\partial_{x_2} f_l\|_{\LL^\infty(\R^2)}\leq \frac{1}{\sigma}\sqrt{\frac{\pi}{2}}$, and by \eqref{YMalDer2} and \eqref{SteinOp} we get Theorem \ref{AsyInd}.  
\medskip

\section{$\LL^2$ Estimates for the KPZ Fixed Point}\label{L2}

In order to apply Malliavin calculus to the KPZ fixed point we need to prove some estimates that will allow us to ensure that the target variables are square integrable. The underlying probability space $\left(\Omega,\cF,\P\right)$ that we work with can be constructed as a product space, $\Omega:=\Omega_1\times\Omega_2$ and $\cF:=\cF_1\otimes\cF_2$, endowed with the product measure $\P:=\P_1\otimes\P_2$, where for every sample $(\fb,\cL)$ we have that $\fb\in\Omega_1$ is a standard two-sided Brownian motion and $\cL\in\Omega_2$ is an independent directed landscape. For the next lemmas, it is worth to recall the definition \eqref{EvoDef} of $\fh_t(x)$ and the definition \eqref{ArgMax} of $Z_t(x)$, where in both cases the initial profile is $\fh_0=\beta\fb$. We start by proving some symmetries of $Z_t(x)$ and by ensuring that it belongs to $\LL^2(\Omega)$.   
 
\begin{lem}\label{ArgSc}
Let $Z$ and $Z_t(x)\equiv Z_{0,t}(x;\beta\fb)$ be given by \eqref{ArgEq} and \eqref{ArgMax}, respectively.
\begin{enumerate}
\item $Z$ is a symmetric random variable;
\item $\P\left[Z=z_0\right]=0$ for all $z_0\in\R$;
\item If $t>0$ and $x\in\R$ are fixed then $Z_t(x)\stackrel{dist.}{=}x+t^{2/3}Z$;
\item $\int_0^\infty u\P\left[|Z|>u\right]du <\infty$.
\end{enumerate}
\end{lem}

\noindent\paragraph{\bf Proof Lemma \ref{ArgSc}}
The symmetry of $Z$ follows from the invariance of the Airy$_2$ process and the two-sided Brownian motion under time reversal $z\mapsto -z$. Now assume that $\P\left[Z=z_0\right]>0$ for some $z_0\in\R$. Let $a,b\in\R$ such that $z_0\in(a,b)$ and 
$$Z_{a,b}:=\argmax_{z\in[a,b]}\left\{\beta\fb(z)+\cA(z)-z^2\right\}=\argmax_{z\in[a,b]}\left\{\beta\left(\fb(z)-\fb(a)\right)+\left(\cA(z)-\cA(a)\right)-z^2\right\}\,.$$
Since $Z$ is the location of the global maximum, $Z_{a,b}=Z$ on the event $Z=z_0\in (a,b)$, and hence 
$$0<\P\left[Z=z_0\right]=\P\left[Z=z_0,Z_{a,b}=Z\right]\leq \P\left[Z_{a,b}=z_0\right]\,.$$   
However, the process $\left(\cA_2(x)-\cA_2(a)\,,\,x\in[a,b]\right)$ is absolutely continuous with respect to Brownian motion \cite{CoHa}, and is independent of $\left(\fb(x)-\fb(a)\,,\,x\in[a,b]\right)$. The location of the maximum of a sum of two independent Brownian motions minus a parabola has a continuous distribution \cite{GrWe}, and hence $\P\left[Z_{a,b}=z_0\right]=0$, which leads to a contradiction. 
\medskip

By \eqref{123} (take $\gamma=t^{-1/3}$), and translation invariance of the Airy sheet,
\begin{eqnarray*}
Z_t(x)&\stackrel{dist.}{=}&\argmax_{z\in\R}\left\{\beta\fb(z)+t^{1/3}\cA((z-x)t^{-2/3})- \frac{(z-x)^2}{t}\right\}\\
&=&\argmax_{z\in\R}\left\{t^{-1/3}\beta\fb(z)+\cA((z-x)t^{-2/3})-((z-x)t^{-2/3})^2\right\}\,.
\end{eqnarray*}
Let $y:=(z-x)t^{-2/3}$, then $Z_t(x)=x+Y_t(x)t^{2/3}$ where
\begin{eqnarray*}
Y_t(x)&=&\argmax_{y\in\R}\left\{t^{-1/3}\beta\fb(x+yt^{2/3})+\cA(y)-y^2\right\}\\
&=&\argmax_{y\in\R}\left\{\beta t^{-1/3}\left(\fb(x+yt^{2/3})-\fb(x)\right)+\cA(y)-y^2\right\}\\
&\stackrel{dist.}{=}&\argmax_{y\in\R}\left\{\beta\fb(y)+\cA(y)-y^2\right\}\,.
\end{eqnarray*}
In the last step we use that for all $x\in\R$ and $t>0$, $t^{-1/3}\left(\fb(x+yt^{2/3})-\fb(x)\right)\stackrel{dist.}{=}\fb(y)$, as process in $y\in\R$, by scaling and shift invariance of the Brownian motion. 
\medskip

Notice that if $u\geq 0$ and $|Z|>u$ then 
$$\cA(0)\leq \max_{z\in\R}\{\beta \fb(z)+\cA(z)-z^2\}= \max_{|z|>u}\{\beta \fb(z)+\cA(z)-z^2\}\,,$$
and hence,
\begin{eqnarray*}
\P\left[|Z|>u\right]&\leq &\P\left[\cA(0)\leq \max_{|z|>u}\{\beta \fb(z)+\cA(z)-z^2\}\right]  \\
&\leq & \P\left[\cA(0)\leq -\frac{u^2}{4}\right]+\P\left[\max_{|z|>u}\{\beta \fb(z)+\cA(z)-z^2\}\geq -\frac{u^2}{4}\right]\,.
\end{eqnarray*}
On one hand, the random variable $\cA(0)$ has a GUE Tracy-Widom distribution and therefore
$$\int_0^\infty u\P\left[\cA(0)\leq -\frac{u^2}{4}\right]<\infty\,.$$
On the other hand, if $\beta \fb(z)\leq z^2/2$ for all $|z|> u$ then $\beta \fb(z)-z^2\leq -z^2/2$ for all $|z|> u$, and hence
\begin{eqnarray*}
\P\left[\max_{|z|>u}\{\beta \fb(z)+\cA(z)-z^2\}\geq -\frac{u^2}{4}\right]&\leq &\P\left[\max_{|z|>u}\{\beta\fb(z)-z^2/2\}>0\right]\\
&+&\P\left[\max_{|z|>u}\{\cA(z)-z^2/2\}\geq -\frac{u^2}{4}\right]\,.
\end{eqnarray*}
Both terms decay to $0$ as $u$ gets large in such a way that one can conclude the proof of property 4 of Lemma \ref{ArgSc}. See for instance Proposition 2.13 \cite{CoLiWa} for the decay of probability involving the Airy$_2$ process minus a parabola, while for the decay of the probability involving the Brownian motion minus a parabola it is a standard estimate.     

\hfill$\Box$\\

\begin{lem}\label{L2KPZ}
For all $a\leq b$ we have that 
\begin{equation}\label{suptail}
\int_0^\infty u\P\left[\sup_{x\in[a,b]}|\fh_t(x)|>u\right]du<\infty\,.
\end{equation}
In particular $\sup_{x\in[a,b]}\fh_t(x)\in\LL^2(\Omega)$.
\end{lem}

\noindent\paragraph{\bf Proof Lemma \ref{L2KPZ}}
For a moment, let us keep track of the dependence on the parameter $\beta>0$ and write $\fh_t\equiv\fh_t^\beta$ and $Z_t\equiv Z_t^\beta$. It is not hard to see that 
\begin{equation*}
\int_0^\infty u\P\left[\sup_{x\in[a,b]}|\fh^{\sqrt{2}}_t(x)|>u\right]du<\infty\,,   
\end{equation*}
since $\fh^{\sqrt{2}}_t(x)=\fh^{\sqrt{2}}_t(0)+\left(\fh^{\sqrt{2}}_t(x)-\fh^{\sqrt{2}}_t(0)\right)$, $\left(\fh^{\sqrt{2}}_t(x)-\fh^{\sqrt{2}}_t(0)\right)$ is a two-sided Brownian motion (time stationarity), and $\fh^{\sqrt{2}}_t(0)\in\LL^2(\Omega)$ (since it has the Baik-Rains distribution). Thus, \eqref{suptail} follows as soon as we prove that 
\begin{equation}\label{SupComp}
\int_0^\infty u\P\left[\sup_{x\in[a,b]}|\fh^{\sqrt{2}}_t(x)-\fh^{\beta}_t(x)|>u\right]du<\infty\,.   
\end{equation}
For $u>0$ define $\fh_t^{\beta,u}(x):=\max_{z\in [-u,u]}\left\{\beta\fb(z)+\cL(z,0;x,t)\right\}$. Let 
$$E_1(u):=\left\{\fh_t^{\beta,u}(x)=\fh_t^{\beta}(x)\,\mbox{ and }\,\fh_t^{\sqrt{2},u}(x)=\fh_t^{\sqrt{2}}(x)\,\, \forall\,\,x\in[a,b]\right\}$$ 
and 
$$E_2(u):=\left\{\sqrt{2}\fb(z)-u/2\leq \beta\fb(z)\leq \sqrt{2}\fb(z)+u/2\,\, \forall\,\,z\in[-u,u]\right\}\,.$$ 
Then 
\begin{equation}\label{SupComp1}
\int_0^\infty u\P\left[E^c_i(u)\right]du<\infty\,,\,\mbox{ for }\,i=1,2\,.
\end{equation}

Indeed, for $E_1(u)$, one has to notice that if $Z^\beta_t(x)\in [-u,u]$ then $\fh^\beta_t(x)=\fh^{\beta,u}_t(x)$, and to use the ordering of the locations (Proposition 2 \cite{Pi1}),
$$Z^\beta_t(a)\leq Z^\beta_t(x)\leq Z^\beta_t(b)\,\,\forall\,x\in [a,b]\,,$$
to see that if $Z^\beta_t(a),Z^\beta_t(b)\in [-u,u]$ then $Z^\beta_t(x)\in [-u,u]$ for all $x\in[a,b]$. Thus
$$\P\left(E^c _1(u)\right)\leq \P\left(|Z^\beta_t(a)|>u\right)+\P\left(|Z^\beta_t(b)|>u\right)+\P\left(|Z^{\sqrt{2}}_t(a)|>u\right)+\P\left(|Z^{\sqrt{2}}_t(b)|>u\right)\,,$$
and by Lemma \ref{ArgSc}, this implies that $\int_0^\infty u\P\left[E^c_1(u)\right]du<\infty$. For $E_2(u)$, one only needs to use classical bounds for the running maximum of a Brownian motion. 
\medskip

To finish the proof, notice that, on the event $E_2(u)$, 
$$\sqrt{2}\fb(z)+\cL(z,0;x,t)-u/2\leq \beta\fb(z)+\cL(z,0;x,t)\leq \sqrt{2}\fb(z)+\cL(z,0;x,t)+u/2\,$$
for all $z\in[-u,u]$, and thus,
$$|\fh^{\sqrt{2},u}_t(x)-\fh^{\beta,u}_t(x)|\leq u\,,$$
This shows that, on the event $E_1(u)\cap E_2(u)$,
$$\sup_{x\in[a,b]}|\fh^{\sqrt{2}}_t(x)-\fh^{\beta}_t(x)|\leq u\,,$$
and therefore, \eqref{SupComp1} implies \eqref{SupComp}. 

\hfill$\Box$\\

\section{Malliavin calculus basics}\label{MalCal}

We work with the isonormal Gaussian process $\left\{\,W(\phi)\,:\,\phi\in\LL^2(\R)\,\right\}$ associated with the standard two-sided Brownian motion on the probability space $(\Omega_1,\cF_1,\P_1)$. In this case  
$$W(\phi):=\int_\R \phi d\fb\,,$$
is the Wiener integral of $\phi$ with respect to $\left(\fb(u)\,,\,u\in\R\right)$. We say that $X\equiv X(\fb)$ is a smooth  random variable if $X\equiv X(\fb)=f\left(W(\phi_1),\cdots,W(\phi_n)\right)$ where $\phi_i\in\LL^2(\R)$ for all $i=1,\dots,n$, and $f:\R^n\to\R$ is a smooth function for which all derivatives have polynomial growth. The Malliavin derivative of a smooth random variable $X$ with respect to $\fb$ is defined as the $\LL^2(\R\times\Omega_1)$ valued random element   
$$u\in\R\,\mapsto\, DX(u):=\sum_{i=1}^n\partial_{x_i} f(W(\phi_i))\phi_i(u)\,.$$
For a smooth random variable $X$ define   
$$\|X\|_{1,p}:=\left(\E_1\left[X^p\right]+\E_1\left[\|DX^p\|_{\LL^2(\R)}\right] \right)^{1/p}\,.$$
For any $p\geq 1$ the derivative operator is closable and its domain can be extended to $\D^{1,p}$, the completion of the space of smooth random variables with respect to $\|\cdot\|_{1,p}$. From now on we take $p=2$, and $\D^{1,2}$ is a Hilbert space with scalar product 
$$\langle X,Y\rangle_{1,2}:=\E_1\left[XY\right]+\E_1\left[\langle DX,DY\rangle_{\LL^2(\R)}\right]\,.$$ 
\medskip

In this paper we use a simplified version of the Malliavin integration by parts formula as follows:
$$\E_1\left[W(\phi)X\right]=\E_1\left[\langle DX,\phi\rangle_{\LL^2(\R)}\right]\,.$$
There is a more general version of the integration by parts formula involving the divergent operator, that is the adjoint of the Malliavin derivative. However, since we are only considering observables that are given by the Wiener integral of a deterministic function, there is no need to introduce the divergent operator. The Malliavin derivative satisfies the following chain rule: if $f:\R^n\to\R$ is a continuously differentiable real valued function with bounded derivatives then 
\begin{equation}\label{IPCR}
D f(X_1,\dots,X_n)= \sum_{k=1}^n\partial_{x_k} f(X_1,\dots,X_n)DX_k\,.
 \end{equation}
If $X\equiv X(\fb,\cL)$ is a square integrable $\cF$-measurable random variable then  
$$\E\left[W(\phi)X\right]=\E\left[\E\left[W(\phi)X\mid \cL\right]\right]\,.$$
By the independence between $\fb$ and $\cL$, and the substitution rule for conditional expectation,
$$\E\left[W(\phi)X\mid \cL=\ff\right]=\E_1\left[W(\phi)X^\ff\right]=\E_1\left[\langle DX^\ff,\phi\rangle_{\LL^2(\R)}\right]\,,$$
where $X^\ff(\cdot)=X(\cdot,\ff):\Omega_1\to\R$. Hence 
\begin{equation}\label{IP}
 \E\left[W(\phi)X\right]=\E\left[\langle DX,\phi\rangle_{\LL^2(\R)} \right]\,.
 \end{equation}
\medskip

Next we are going to give some examples and use the following function several times:
\begin{equation}\label{zeta}
u\in\R\,\mapsto\,\zeta_x(u)=\left\{\begin{array}{ll}\,\,\,\,\1_{(0,x]}(u)& \mbox{ if } x>0\,,\\
\,\,\,\,0& \mbox{ if } x=0\,,\\
- \1_{(x,0]}(u) &\mbox{ if } x<0\,.\end{array}\right.
\end{equation}
A simple computation shows that if $x\leq y$ then 
$$\zeta_y(u)-\zeta_x(u)=\1_{(x,y]}(u)\,. $$
By definition, we have that 
\begin{equation}\label{ExMalDer1}
\mbox{ if }X=\int_\R\phi d(\beta\fb)+c=\beta W(\phi)+c\,\,\,\mbox{ then }\,\,\,DX(u)=\beta\phi(u)\,
\end{equation}  
(take $f(x)=\beta x+c$). Since $\fb(x)=W\left(\zeta_x\right)$,  
\begin{equation}\label{ExMalDer2}
\mbox{ if }X=\beta\fb(x)+c\mbox{ then }DX(u)=\beta\zeta_x(u)\,.
\end{equation}  
Another key example in our context is the Malliavin derivative of the maximum as follows. Assume that $\ff:\R\to\R$ is a continuous function such that the maximum of 
$$\fb_\ff(z):=\beta\fb(z)+\ff(z)\,,\,\mbox{ for $z\in\R$}\,,\,$$
is in $\LL^2(\Omega_1)$ and it is attained  $\P_1$-a.s at a unique location $\tau=\argmax_{z\in\R}\left\{\fb_\ff(z)\right\}$. If $\E_1\left[|\tau|\right]<\infty$ then (recall the definition \eqref{zeta} of $\zeta_x$)
\begin{equation}\label{ExMalDer3}
M=\max_{z\in\R}\left\{\fb_\ff(z)\right\}\in\D^{1,2}\,\,\mbox{ and }\,\,DM(u)=\zeta_\tau(u)\,.
\end{equation}
To justify \eqref{ExMalDer3}, we assume without loss of generality that $\ff(0)=0$. Let $a>0$ and consider the maximum $M^a$ of $\fb_\ff(z)$ over $z\in[-a,a]$, and assume that it is attained  $\P_1$-a.s at a unique location $\tau^a$.  We can approximate $M^a$ by the maximum over a finite set such that $0\in\{z_1,\dots,z_n\}\nearrow\Q\cap [-a,a]$:
$$M^a_n:=\max_{k=1,\dots,n}\left\{X_k\right\}\mbox{ and }\tau^a_n:=\argmax_{k=1,\dots,n}\left\{\fb_\ff(z_k)\right\}\,,$$
where $X_k:=\fb_\ff(z_k)\in\D^{1,2}$ for each $k=1,\dots,n$, by \eqref{ExMalDer2} (take $x=z_i$ and $c=\ff(z_i)$). The function $f_n(x_1,\dots,x_n):=\max_{k=1,\dots,n}\left\{x_i\right\}$ is not continuously differentiable, but it is Lipschitz continuous and its partial derivatives exist almost everywhere (with respect to the Lebesgue measure on $\R^n$), which allow us to apply smoothing arguments to use the chain rule \eqref{IPCR} (Proposition 4.2 \cite{Nu}). Thus $M^a_n=f_n\left(X_1,\dots,X_n \right)\in\D^{1,2}$  and 
$$DM^a_n(u)=\sum_{k=1}^n\partial_{x_k}f_n\left(X_1,\dots,X_n \right)DX_k(u)\,.$$ 
Let $A_1:=\left\{f_n(x_1,\dots,x_n)=x_1\right\}$ and for $k=2,\dots,n$ let
$$A_k=\left\{f_n(x_1,\dots,x_n)\neq x_1,\dots,f_n(x_1,\dots,x_n)\neq x_{k-1}\,,\,f_n(x_1,\dots,x_n)=x_k\right\}\,.$$
Then $\partial_{x_k}f_n(x_1,\dots,x_n)=\1_{A_k}(x_1,\dots,x_n)$ almost everywhere, and together with \eqref{ExMalDer2}, this implies that
$$DM^a_n(u)=\sum_{k=1}^n \1_{A_k}\left(X_1,\dots,X_n \right) \zeta_{z_k}(u)=\zeta_{\tau_n^a}(u).$$
Now we use that $\fb_\ff(0)=0$ and hence $0\leq M^a_n\leq M^a$ and $0\leq M^a-M^a_n\leq M^a\leq M$. By continuity of $\fb_\ff$,  $M^a_n\to M^a$ $\P_1$-a.s. and, by dominated convergence, $M^a_n\to M^a$ in $\LL^2(\Omega_1)$. Since
$$\int_\R|DM^a_n-\zeta_{\tau^a}|^2 du=\int_\R|\zeta_{\tau^a_n}-\zeta_{\tau^a}| du=|\tau^a_n-\tau^a |\leq 2a\,,$$  
by continuity of $\fb_\ff$, $\tau_n^a\to\tau^a$ $\P_1$-a.s., and by dominated convergence, $\tau_n^a\to\tau^a$ in $\LL^1(\Omega_1)$. Hence
$$\E_1\left[\int_\R|DM^a_n-\zeta_{\tau^a}|^2du \right]=\E_1\left[|\tau^a_n-\tau^a |\right]\to 0\,,$$
and we can conclude that $M^a\in\D^{1,2}$ and that $DM^a=\zeta_{\tau^a}$. Now, since $M^a\to M$ and $\tau^a\to\tau$ $\P_1$-a.s., $0\leq M-M^a\leq M$ and (notice that $\tau^a=\tau$ if $|\tau|\leq a$, and that $|\tau^a-\tau|\leq 2\tau$ if $|\tau|>a$) 
$$\int_\R|DM^a-\1_{\tau}|^2 du=\int_\R|\zeta_{\tau^a}-\zeta_{\tau}| du=|\tau^a-\tau |\1_{|\tau|>a}\leq 2|\tau|\1_{\{|\tau|>a\}}\leq 2|\tau|\,,$$ 
we can use dominated convergence again (recall that $\E_1\left[|\tau|\right]<\infty$), to conclude the proof of \eqref{ExMalDer3}.
\medskip

For the nex lemmas it is worth to recall the definitions \eqref{EvoDef}, \eqref{ArgMax} and \eqref{DefObs} of $\fh_t(x)$, $Z_t(x)$ and $X_t^\phi$, respectevely.
\begin{lem}\label{BasMalDer}
We have that $X_0^\phi\in\D^{1,2}$ and $DX_0^\phi(u)=\beta\phi(u)$. Furthermore, for $t>0$ we have that $\fh_t(x)\in\D^{1,2}$ $\P_2$-a.s. and $D(\fh_t(x))(u)=\beta\zeta_{Z_t(x)}(u)$.
\end{lem}

\noindent\paragraph{\bf Proof Lemma \ref{BasMalDer}}
By Lemma \ref{ArgSc} and Lemma \ref{L2KPZ}, $\E\left[\fh_t(x)^2\right]<\infty$ and $\E\left[|Z_t(x)|\right]<\infty$. Thus, $\P_2$-a.s.  $\E_1\left[\fh_t(x)^2\right]<\infty$ and $\E_1\left[|Z_t(x)|\right]<\infty$. Since
$$X_0^\phi=W(\beta\phi)=\beta\int_\R\phi d\fb\,\mbox{ and }\,\fh_t(x)=\max_{z\in\R}\left\{\beta\fb(z)+\cL(z,0;x,t)\right\}\,,$$ 
we have that Lemma \ref{BasMalDer} follows from \eqref{ExMalDer1} and \eqref{ExMalDer3}. 

\hfill$\Box$\\

\begin{lem}\label{ConvFor}
Let $\phi\in\SS(\R)$ and $t>0$. We have that $X_t^{\phi}\in\D^{1,2}$ $\P_2$-a.s. and 
$$\|DX_t^{\phi}\|_{\LL^2(\Omega\times \R)}=\beta\|\phi\|_{\LL^2(\R)}\,.$$
Furthermore, for all $\phi_1\in\LL^2(\R)$ and $\phi_2\in\SS(\R)$ we have that 
$$\E\left[\langle DX_t^{\phi_2},\phi_1\rangle_{\LL^2(\R)}\right]=\beta\E\left[(\phi_1\star\phi_2)(t^{2/3}Z)\right]\,.$$
\end{lem}

\noindent\paragraph{\bf Proof Lemma \ref{ConvFor}} 
Denote $\Delta_xf(y):=f(y)-f(x)$ and $\phi(x)=\sum_{j=1}^n c_j \1_{(x_{j-1},x_j]}(x)$. Since $Z_t(x)$ is a nondecreasing function of $x\in\R$  (Proposition 2 \cite{Pi1}), by Lemma \ref{BasMalDer},
$$D\left(\Delta_x\fh_t(y)\right)(u)=D\left(\fh_t(y)\right)(u)-D\left(\fh_t(x)\right)(u)=\beta\zeta_{Z_t(y)}(u)-\beta\zeta_{Z_t(x)}(u)=\beta \1_{\left(Z_t(x),Z_t(y)\right]}(u)\,,$$
for $x<y$. Therefore, if $\phi\in\SS(\R)$ then  $X_t^{\phi}:=\int_\R\phi \partial_x\fh_t:=\sum_{j=1}^n c_j \Delta_{x_{j-1}}\fh_t(x_j)\in\D^{1,2}$ and 
$$DX_t^{\phi}(u)=\sum_{j=1}^n c_j D\left(\Delta_{x_{j-1}}\fh_t(x_j)\right)(u)=\beta\sum_{j=1}^n c_j \1_{ \left(Z_t(x_{j-1}),Z_t(x_{j})\right]}(u)\,.$$
Thus, $\left(DX_t^{\phi}(u) \right)^2=\beta^2\sum_{j=1}^n c_j^2\1_{ \left(Z_t(x_{j-1}),Z_t(x_{j})\right]}(u)$ and, by Lemma \ref{ArgSc},
\begin{eqnarray*}
\|DX_t^{\phi}\|^2_{\LL^2(\Omega\times \R)}&=&\E\left[\int_\R\left(DX_t^{\phi}(u) \right)^2du\right]\\
&=&\beta^2\sum_{j=1}^n c_{j}^2\left(\E \left[Z_t(x_j)\right]-\E\left[Z_t(x_{j-1})\right]\right)\\
&=&\beta^2\sum_{j=1}^n c_{j}^2\left(x_j-x_{j-1}\right)\\
&=& \beta^2\|\phi\|^2_{\LL^2(\R)}\,.
\end{eqnarray*}
Denote $\psi_1(x):=\int_0^x\phi_1(u)du$ and $\phi_2(x)=\sum_{j=1}^n c_j \1_{(x_{j-1},x_j]}(x)$. Thus,
\begin{eqnarray*}
\langle DX_t^{\phi_2},\phi_1\rangle_{\LL^2(\R)}&=&\int_\R DX_t^{\phi_2}(u)\phi_1(u)du\\
&=&\beta\sum_{j=1}^nc_{j} \int_{Z_t(x_{j-1})}^{Z_t(x_{j})}\phi_1(u) du\\
&=&\beta\sum_{j=1}^n c_{j} \left(\psi_1(Z_t(x_j))-\psi_1(Z_t(x_{j-1}))\right)\,.
\end{eqnarray*}
By Lemma \ref{ArgSc},
\begin{eqnarray*}
\E\left[\langle DX_t^{\phi_2},\phi_1\rangle_{\LL^2(\R)}\right]&=&\beta\sum_{j=1}^n c_{j} \Big(\E\left[\psi_1(Z_t(x_j))\right]-\E\left[\psi_1(Z_t(x_{j-1}))\right]\Big)\\ 
&=&\beta\sum_{j=1}^n c_{j} \Big(\E\left[\psi_1(x_j+t^{2/3}Z)\right]-\E\left[\psi_1(x_{j-1}+t^{2/3}Z)\right]\Big)\\
&=&\beta \E\left[\sum_{j=1}^nc_{j} \Big(\psi_1(x_j+t^{2/3}Z)-\psi_1(x_{j-1}+t^{2/3}Z)\Big)\right]\\
&=& \beta\E\left[\sum_{j=1}^n c_{j} \int_{x_{j-1}+t^{2/3}Z}^{x_j+t^{2/3}Z}\phi_1(u) du\right]\\
&=& \beta\E\left[\int_\R\phi_1(u)\left(\sum_{j=1}^n c_{j} \1_{\left(x_{j-1}+t^{2/3}Z,x_j+t^{2/3}Z\right]}(u)\right) du\right]\\
&=& \beta\E\left[\int_\R\phi_1(u)\phi_2(u-t^{2/3}Z) du\right]\\
&=& \beta\E\left[\int_\R\phi_1(u)\phi_2(u+t^{2/3}Z)du\right]\,,
\end{eqnarray*}
which shows that 
$$ \E\left[\langle DX_t^{\phi_2},\phi_1\rangle_{\LL^2(\R)}\right]=\beta\E\left[\phi_1\star\phi_2(t^{2/3}Z)\right]\,,$$
according to the definition \eqref{CrossCorr} of the cross correlation of $\phi_1$ and $\phi_2$.

\hfill$\Box$\\
 
Now we consider the cases where $\beta\neq \sqrt{2}$ and $\phi\in\CC^1_b(\R)$, or $\beta=\sqrt{2}$ and $\phi\in\LL^2(\R)$. The operator $D$ is closable, and to get the existence of $DX_t^{\phi}$ we show that there exists a sequence $X_{n}\in\D^{1,2}$ such that $\P_2$-a.s. $X_{n}$  converges to $X_t^{\phi}$ in $\LL^2(\Omega_1)$ and $DX_{n}$ converges in $\LL^2(\Omega_1\times\R)$ to some element $G$. In this case, $X_t^{\phi}\in \D^{1,2}$ and $DX_t^{\phi}=G$. 

\begin{lem}\label{Approx}
Consider the cases: (i) $\beta\neq \sqrt{2}$ and $\phi\in\CC^1_b(\R)$; (ii) $\beta=\sqrt{2}$ and $\phi\in\LL^2(\R)$. For all $t>0$ there exists a sequence $\phi_n\in\SS(\R)$ for $n\geq 1$ such that
\begin{equation}\label{iso1}
\lim_{n\to\infty}\|\phi_n-\phi\|_{\LL^2(\R)}=0\,\mbox{ and }\,\lim_{n\to\infty}\|X_t^{\phi_n}-X_t^{\phi}\|_{\LL^2(\Omega)}=0\,.
\end{equation}
\end{lem}

\noindent\paragraph{\bf Proof Lemma \ref{Approx}}
For $\beta=\sqrt{2}$ and $\phi\in\LL^2(\R)$, take a sequence $\phi_n\in\SS(\R)$   such that
$$\lim_{n\to\infty}\|\phi_n-\phi\|_{\LL^2(\R)}=0\,.$$ 
By isometry, 
\begin{equation*}
\lim_{n\to\infty}\|X_t^{\phi_n}-X_t^{\phi}\|_{\LL^2(\Omega)}=\sqrt{2}\lim_{n\to\infty}\|\phi_n-\phi\|_{\LL^2(\R)}=0\,.
\end{equation*}
For $\beta\neq \sqrt{2}$ and $\phi\in\CC^1_b(\R)$, consider real numbers $a<b$ such that $\supp(\phi)\subseteq (a,b)$, and let  $x^n_0=x^n_1=a<x^n_2<\dots<x^n_{n-1}<x^n_n=b$ be a partition of $[a,b]$, such that $\max_{j=1,\dots,n}(x^n_{j}-x^n_{j-1})\to 0$ as $n\to \infty$, and define $\phi_n(z):=\sum_{j=1}^{n}\phi(x^n_j)\1_{(x^n_{j-1},x^n_{j}]}(z)$. Thus, we clearly have that $\|\phi_n-\phi\|_{\LL^2(\R)}\to 0$ as $n\to\infty$. On the other hand, 
$$X_t^{\phi}:=-\int_{\R}\phi'(z)\fh_t(z)dz=-\int_{[a,b]}\phi'(z)\fh_t(z)dz\,,$$
and by the mean value theorem, there exists $z^n_j\in(x^n_{j},x^n_{j+1}]$ for each $j=1,\dots,n-1$ such that   
$$X_t^{\phi_n}=\sum_{j=1}^{n}\phi(x^n_j)\left(\fh_t(x^n_{j})-\fh_t(x^n_{j-1})\right)=-\sum_{j=1}^{n-1}\fh_t(x^n_{j})\phi'(z^n_j)(x^n_{j+1}-x^n_{j})\,.$$
By continuity of $\fh_t$ and $\phi'$, $X_t^{\phi_n}$ is converging a.s. to $X_t^{\phi}$, as $n\to\infty$. Furthermore, 
$$\big|X_t^{\phi_n}\big|=\big|\sum_{j=1}^{n-1}\fh_t(x^n_{j})\phi'(z^n_j)(x^n_{j+1}-x^n_{j})\big|\leq C\sup_{x\in[a,b]}|\fh_t(x)|\,,$$
where $C=\sup_{x\in [a,b]}|\phi'(x)|\left(b-a\right)$. By Lemma \ref{L2KPZ}, $\sup_{x\in[a,b]}|\fh_t(x)|\in\LL^2(\Omega)$, and by dominated convergence, $\|X_t^{\phi_n}-X_t^{\phi}\|_{\LL^2(\Omega)}\to 0$ as $n\to\infty$.  

\hfill$\Box$\\

\begin{lem}\label{ExiMalDer}
Consider the cases: (i) $\beta\neq \sqrt{2}$ and $\phi\in\CC^1_b(\R)$; (ii) $\beta=\sqrt{2}$ and $\phi\in\LL^2(\R)$. Then for all $t>0$ we have that $X_t^{\phi}\in\D^{1,2}$ $\P_2$-a.s. and 
$$\|DX_t^{\phi}\|_{\LL^2(\Omega\times \R)}=\beta\|\phi\|_{\LL^2(\R)}\,.$$
Furthermore, for $\phi_1\in\LL^2(\R)$ and $\phi_2$ as in (i) or (ii), we have that 
$$\E\left[\langle DX_t^{\phi_2},\phi_1\rangle_{\LL^2(\R)}\right]=\beta\E\left[(\phi_1\star\phi_2)(t^{2/3}Z)\right]\,.$$
\end{lem}

\noindent\paragraph{\bf Proof Lemma \ref{ExiMalDer}}
Let $\phi_n$ for $n\geq 1$ be a sequence in $\SS(\R)$, as in Lemma \ref{Approx}. Given $\phi_n$ and $\phi_m$  chose a refinement of both partitions to write 
$$\phi_n(z)=\sum_{j=1}^k c_{n,j} \1_{(x_{j-1},x_{j}]}(z)\,\,\mbox{ and }\,\,\phi_m(z)=\sum_{j=1}^k c_{m,j} \1_{(x_{j-1},x_{j}]}(z)\,.$$
Then
$$DX_t^{\phi_n}(z)=\beta\sum_{j=1}^k c_{n,j} \1_{ \left(Z_t(x_{j-1}),Z_t(x_{j})\right]}(z)\,\,\mbox{ and }\,\,DX_t^{\phi_m}(z)=\beta\sum_{j=1}^k c_{m,j} \1_{ \left(Z_t(x_{j-1}),Z_t(x_{j})\right]}(z)\,.$$
By Lemma \ref{ConvFor}, 
\begin{equation}\label{MallIso}
\|DX_t^{\phi_n}-DX_t^{\phi_m}\|_{\LL^2(\Omega\times \R)}=\|DX_t^{\phi_n-\phi_m}\|_{\LL^2(\Omega\times \R)}=\beta\|\phi_n-\phi_m\|_{\LL^2(\R)}\,.
\end{equation}
Since $\phi_n$ converges to $\phi$ in $\LL^2(\R)$, \eqref{MallIso} implies that $DX^{\phi_n}_{t}$ is a Cauchy sequence in $\LL^2(\Omega\times \R)$. Therefore, there exists $G=(G(x)\,;\,x\in\R)\in\LL^2(\Omega\times \R)$ such that 
\begin{equation}\label{iso2}
\lim_{n\to\infty}\|DX_t^{\phi_n}-G\|_{\LL^2(\Omega\times \R)} =0\,.
\end{equation}
By \eqref{iso1}, \eqref{iso2} and Fubini's theorem, we can conclude that there exists a subsequence $\phi_{n_k}\in\SS(\R)$ for $k\geq 1$ such that $\P_2$-a.s. 
$$\lim_{k\to\infty}\|X_t^{\phi_{n_k}}-X_t^{\phi}\|_{\LL^2(\Omega_1)}=0\,,$$
and 
$$\lim_{k\to\infty}\|DX_t^{\phi_{n_k}}-G\|_{\LL^2(\Omega_1\times \R)} =0\,.$$
Thus, $\P_2$-a.s. $X_t^{\phi}\in \D^{1,2}$ and $DX_t^{\phi}=G$. By approximating $\phi$ with $\phi_n\in\SS(\R)$, and applying Lemma \ref{ConvFor}, we have that   
$$\|DX_t^{\phi}\|_{\LL^2(\Omega\times \R)}=\beta\|\phi\|^2_{\LL^2(\R)}\,.$$
Similarly (by approximation), 
$$\E\left[\langle DX_t^{\phi_2},\phi_1\rangle_{\LL^2(\R)}\right]=\beta\E\left[\left(\phi_1\star\phi_2\right)(t^{2/3}Z)\right]\,,$$
for all $\phi_1\in\LL^2(\R)$ and $\phi_2$ as in (i) or (ii).

\hfill$\Box$\\

\begin{rem}\label{Correction}
Another way to prove Lemma \ref{ExiMalDer} is to notice that (Lemma \ref{BasMalDer}) if $\phi_2\in\CC^1_b(\R)$ then 
$$DX_t^{\phi_2}(u)=-\beta\int_{\R}\phi_2'(z)\zeta_{Z_t(z)}(u)dz \,.$$
By Fubini's Theorem, this implies that 
\begin{eqnarray*}
\nonumber\langle DX_t^{\phi_2},\phi_1\rangle_{\LL^2(\R)}&=&\int_\R DX_t^{\phi_2}(u)\phi_1(u)du\\
\nonumber&=&-\beta\int_\R\left(\int_\R\phi'_2(z)\zeta_{Z_t(z)} (u)\phi_1(u)dz\right)du\\
&=&-\beta\int_\R\phi'_2(z)\psi_1 (Z_t(z))dz\,.
\end{eqnarray*}
(Recall that $\psi'_1(z)=\int_\R\phi_1(u)\zeta_z(u)du$.) Together with Lemma \ref{ArgSc} (and standard integration by parts), this implies Lemma \ref{ExiMalDer}. It also implies that 
$$\E\left[\left|\langle DX_t^{\phi_2},\phi_1\rangle_{\LL^2(\R)}\right|\right]\leq \beta\E\left[|\psi_1|\star|\phi'_2|(t^{2/3}Z)\right]\,.$$ 

\end{rem}

\begin{lem}\label{LemVarDen}
Let $g_\beta(x):=\Var\left[\fh_1(x)\right]$ and $\bF_\beta(x):=\P\left[Z\leq x\right]$. Then $g_\beta$ is differentiable and 
$$g'_\beta(x)=\beta^2\left(2\bF_\beta(x)-1\right)\,.$$
\end{lem}

\noindent\paragraph{\bf Proof Lemma \ref{LemVarDen}}
Take $u=z-x$ and use shift invariance of Brownian motion together with space stationarity of the Airy$_2$ process to have that
\begin{eqnarray*}
\fh_1(x)-\fh_0(x)&\stackrel{dist.}{=}&\max_{z\in\R}\left\{\fh_0(z)-\cA_2(z)-(z-x)^2\right\}-\fh_0(x)\\
&=&\max_{u\in\R}\left\{\left(\fh_0(u+x)-\fh_0(x)\right)-\cA_2(u+x)-u^2\right\}\\
&\stackrel{dist.}{=}&\max_{u\in\R}\left\{\fh_0(u)-\cA_2(u)-u^2\right\}\\
&\stackrel{dist.}{=}&\fh_1(0)\,.
\end{eqnarray*}
Notice that this distributional equality holds for all $\beta>0$ and $x\in \R$ fixed, although \eqref{stat_1} holds only for $\beta^2=2$. Write $\fh_1(x)=\left(\fh_1(x)-\fh_0(x)\right)+\fh_0(x)$, and use that $\fh_1(x)-\fh_0(x)\stackrel{dist.}{=}\fh_1(0)$ (and that $\E\left[\fh_0(x)\right]=0$), to show that  
\begin{eqnarray*}
g_\beta(x)&=&\Var\left[\fh_1(x)\right]\\
&=&\Var\left[\fh_1(x)-\fh_0(x)\right]+\Var\left[\fh_0(x)\right]+2\Cov\left[\fh_0(x),\fh_1(x)-\fh_0(x)\right]\\
&=&\Var\left[\fh_1(0)\right]-\Var\left[\fh_0(x)\right]+2\Cov\left[\fh_0(x),\fh_1(x)\right]\\
&=&\Var\left[\fh_1(0)\right]-\beta^2 |x|+2\E\left[\fh_0(x)\fh_1(x)\right]\,.
\end{eqnarray*}

At a first glance, it seems that we may have a problem with the derivative at zero because of the modulus function. However, the covariance between $\fh_0(x)$ and $\fh_1(x)$, together with the symmetry of $Z$, is going to compensate that. Notice that  (recall the definition \eqref{zeta} of $\zeta_x$)
$$\fh_0(x)=\int_{\R}\zeta_x\partial_x\fh_0dx=\beta W(\zeta_x)\,.$$
Hence,  by \eqref{IP} we get that, 
$$\E\left[\fh_1(x)\fh_0(x)\right]=\beta\E\left[\langle D(\fh_1(x)),\zeta_x\rangle_{\LL^2(\R)}\right]\,.$$  
For $y\in\R$ denote $y_+:=\max\{0,y\}$ and $y_-:=\min\{0,y\}$. By Lemma \ref{ArgSc} and Lemma \ref{BasMalDer}, for $x\geq 0$ we have 
$$\E\left[\langle D(\fh_1(x)),\zeta_x\rangle_{\LL^2(\R)}\right]=\beta\E\left[ \psi^+_x(x+Z)\right]\,,$$  
where $\psi^+_x(y)=\int_0^y\zeta_x(z)dz=\min\{x,y_+\}$, while for $x\leq 0$ we have 
$$\E\left[\langle D(\fh_1(x)),\zeta_x\rangle_{\LL^2(\R)}\right]=\beta\E\left[\psi^-_x(x+Z)\right]\,,$$ 
where $\psi^-_x(y)=\int_0^y\zeta_x(z)dz=-\max\{x,y_-\}$ (recall that $\int_0^y=-\int_y^0$ for $y<0$). Since $\psi^+_0(y)=\psi^-_0(y)=0$ for all $y\in\R$, we can write
\begin{equation}\label{EqDer}
g_\beta(x)=\Var\left[\fh_1(0)\right] -\beta^2|x|+2\beta^2 E(x)\,,
\end{equation}
where 
$$E(x):=\left\{\begin{array}{ll}\E\left[ \psi^+_x(x+Z)\right] & \mbox{ if } x\geq 0\,,\\
 \E\left[\psi^-_x(x+Z)\right]&\mbox{ if } x\leq 0\,.\end{array}\right.$$ 
\medskip 

For $x\geq 0$ we get that
$$\psi^+_x(x+Z)= \left\{\begin{array}{lll}  0 &\mbox{ if } Z\leq -x\,,\\
 x+Z &\mbox{ if } Z\in (-x,0]\,,\\
 x &\mbox{ if } Z>0 \,.\end{array}\right.$$ 
Thus,
$$E(x)=\E\left[ \psi^+_x(x+Z)\right]=\E\left[Z\1_{\{Z\in (-x,0]\}}\right]+x\P\left[Z>-x\right]\,,$$
and if $0\leq x\leq y$ then  
\begin{equation}\label{DerEq1}
E(y)-E(x)=\E\left[(y+Z)\1_{\{Z\in (-y,-x]\}}\right]+(y-x)\P\left[Z>-x\right]\,.
\end{equation}
For $Z\in (-y,-x]$ we have that $0\leq y+Z\leq y-x$ and hence
$$0\leq \E\left[(y+Z)\1_{\{Z\in (-y,-x]\}}\right]\leq (y-x)\P\left[ Z\in (-y,-x]\right]\,.$$
By Lemma \ref{ArgSc}, $\P\left[Z=z_0\right]=0$ for all $z_0\in\R$ and \eqref{DerEq1} implies that for all $x>0$
$$E'(x)=\P\left[Z>-x\right]=\P\left[Z\leq x\right]\,,$$
where we also use symmetry of $Z$ in the last equality. Hence, by \eqref{EqDer},
$$g_\beta'(x)=-\beta^2+2\beta^2\P\left[Z\leq x\right]=\beta^2\left(2\bF_\beta(x)-1\right)\,,\mbox{ for $x>0$}\,.$$
For $x\leq 0$ we have that 
$$\psi^-_x(x+Z)= \left\{\begin{array}{lll}  -x &\mbox{ if } Z\leq 0\,,\\
 -(x+Z) &\mbox{ if } Z\in (0,-x]\,,\\
 0 &\mbox{ if } Z> -x \,.\end{array}\right.$$ 
Thus,
$$E(x)=\E\left[\psi^-_x(x+Z)\right]=-x\P\left[Z\leq -x\right]-\E\left[Z\1_{\left\{ Z\in (0,-x]\right\}}\right]\,,$$
and if $x\leq y\leq 0$ then  
\begin{equation}\label{DerEq2}
E(y)-E(x)=\E\left[(y+Z)\1_{\{Z\in (-y,-x]\}}\right]-(y-x)\P\left[Z\leq-x\right]\,.
\end{equation}
By using Lemma \ref{ArgSc} again, and \eqref{DerEq2}, we can deduce that for all $x<0$ 
$$E'(x)=-\P\left[Z\leq -x\right]=\P\left[Z\leq x\right]-1\,,$$ 
and by \eqref{EqDer},
$$g_\beta'(x)=\beta^2+2\beta^2\left(\P\left[Z\leq x\right]-1\right)= \beta^2\left(2\bF_\beta(x)-1\right)\,,\mbox{ for $x<0$}\,.$$

Now we consider $x=0$. By \eqref{DerEq1}, if $z\geq 0$ then (take $x=0$ and $z=y$)
\begin{eqnarray*}
g_\beta(z)-g_\beta(0)&=&-\beta^2 z+2\beta^2\left(E(z)-E(0)\right)\\
&=&-\beta^2z+2\beta^2\left(\E\left[(z+Z)\1_{\{Z\in (-z,0]\}}\right]+z\P\left[Z>0\right]\right)\\
&=& 2\beta^2\E\left[(z+Z)\1_{\{Z\in (-z,0]\}}\right]\,.
\end{eqnarray*}
By \eqref{DerEq2}, if $z\leq 0$ then (take $x=z$ and $y=0$)
\begin{eqnarray*}
g_\beta(z)-g_\beta(0)&=&\beta^2 z+2\beta^2\left(E(z)-E(0)\right)\\
&=&\beta^2z-2\beta^2\left(\E\left[Z\1_{\{Z\in (0,-z]\}}\right]+z\P\left[Z\leq-z\right]\right)\\
&=&\beta^2 z\left(1-2\P\left[Z\leq-z\right]\right)-2\beta^2\E\left[Z\1_{\{Z\in (0,-z]\}}\right]\,.
\end{eqnarray*}
Thus, by Lemma \ref{ArgSc}, 
$$\exists\,\lim_{z\to 0}\frac{g_\beta(z)-g_\beta(0)}{z}=0=\beta^2 \left(2\bF_\beta(0)-1\right)\,.$$

\hfill$\Box$\\

\noindent\paragraph{\bf Proof of Theorem \ref{CovArg}} 
By integration by parts \eqref{IP} we have \eqref{IntPart} and, together with Lemma \ref{ExiMalDer}, this implies the theorem.

\hfill$\Box$\\

\noindent\paragraph{\bf Proof of Theorem \ref{VarDen}} 
It follows directly from Lemma \ref{LemVarDen}, and the covariance can be computed using Theorem \ref{CovArg}. 

\hfill$\Box$\\

\section{Stein's method basics}\label{SteMet}

To contextualize the main idea we start by a brief discussion of the classical Stein's method for normal approximations. Define the functional operator $\cN$, acting on differentiable functions, by 
$$\cN f(x):=\sigma^2f'(x)-xf(x)\,.$$
Recall that  $X\sim N(0,\sigma^2)$ if $X$ has the normal distribution with zero mean and variance $\sigma^2$. The next result is called ``Stein's lemma'' \cite{Ro}:
\begin{itemize}
\item If $X\sim N(0,\sigma^2)$ then $\E\left[\cN f(X)\right]=0$ for all absolute continuous $f$ with $\E\left[f'(X)\right]<\infty$. In particular (integration by parts formula),
$$\E\left[Xf(X)\right]=\sigma^2\E\left[f'(X)\right]\,.$$  
\item If for some random variable $X$, $\E\left[\cN f(X)\right]=0$ for all absolute continuous $f$ with bounded derivative then $X\sim N(0,\sigma^2)$.
\end{itemize}
Motivated by Stein's lemma, one can think that if $\E\left[\cN f(X)\right]$ is close to zero, then $\P_X$ should be close to $N(0,\sigma^2)$. 
\medskip

Turning back to our context, the point is not a normal approximation but asymptotic independence. In this direction, the following two claims can be seen as a two-dimensional version of Stein's lemma, and besides the characterisation of the distribution of $X_1$, it includes independence between the components of a random vector $(X_1,X_2)$. As far as the author's knowledge goes, it has never appeared in the mathematical literature before, and we give a brief explanation as follows. Let $\sigma>0$ and define the functional operator $\cN$, acting on differentiable functions $f:\R^2\to\R$, by
\begin{equation}\label{SteinChar}
\cN f(x_1,x_2):=\sigma^2\partial_{x_1} f(x_1,x_2)- x_1 f(x_1,x_2)\,.
\end{equation}
\begin{itemize}
\item If $X_1\sim N(0,\sigma^2)$ is independent of $X_2$ then 
$$\E\left[\cN f(X_1,X_2)\right]=0\,,$$ 
for all differentiable $f$ with $\E\left[\partial_{x_1} f(X_1,X_2)\right]<\infty$. In particular,
$$\E\left[X_1f(X_1,X_2)\right]=\sigma^2\E\left[\partial_{x_1} f(X_1,X_2)\right]\,.$$  
\item If for some random vector $(X_1,X_2)$, with $\E\left[|X_1|\right]<\infty$, we have that $\E\left[\cN f(X_1,X_2)\right]=0$ for all differentiable $f$ with $\|\partial_{x_1} f\|_{\LL^\infty(\R^2)}<\infty$ then $X_1\sim N(0,\sigma^2)$ and $X_1$ is independent of $X_2$.
\end{itemize}
The first claim follows easily using Fubini's theorem together with (one dimensional) Stein's lemma. For the second claim consider the characteristic function $\psi_{X_1,X_2}$ of $(X_1,X_2)$. Since $\E\left[|X_1|\right]<\infty$ we have 
\begin{eqnarray*}
\partial_{\lambda_1}\psi_{X_1,X_2}(\lambda_1,\lambda_2)&=&i\E\left[ X_1e^{i\left(\lambda_1 X_1+\lambda_2 X_2\right)}\right]\\  
&=&i\sigma^2\E\left[ \partial_{x_1} e^{i\left(\lambda_1 X_1+\lambda_2 X_2\right)}\right]\\
&=&-\lambda_1\sigma^2\psi_{X_1,X_2}(\lambda_1,\lambda_2)\,,
\end{eqnarray*}
where we use in the second equality that $\E\left[\cN f(X_1,X_2)\right]=0$ for $f(x_1,x_2)=\cos\left(\lambda_1 x_1+\lambda_2 x_2\right)$ and for $f(x_1,x_2)=\sin\left(\lambda_1 x_1+\lambda_2 x_2\right)$. Since $\psi_{X_1,X_2}(0,\lambda_2)\equiv\psi_{X_2}(\lambda_2)$ we must have that 
$$\psi_{X_1,X_2}(\lambda_1,\lambda_2)=e^{-\frac{\sigma^2}{2}\lambda_1^2}\psi_{X_2}(\lambda_2)\,,$$
which implies independence between $X_1$ and $X_2$.

Now, to bound the difference between the joint law $\theta\equiv\P_{X_1,X_2}$ and $\eta\equiv\P_{X_1}\otimes\P_{X_2}$, we follow Stein's idea and look at a suitable solution of the partial differential equation
\begin{equation}\label{SteinEq}
\cN f(x_1,x_2)= l(x_1,x_2)-\E\left[l(X_1,x_2)\right]\,,
\end{equation}
where $l:\R^2\to\R$ is a continuously differentiable real valued function with bounded partial derivatives. By Fubini's theorem,
$$\int_{\R^2} \E\left[l(X_1,x_2)\right]d\theta=\int_{\R}\left(\int_\R l(x_1,x_2)d\P_{X_1}(x_1)\right)d\P_{X_2}(x_2)=\int_{\R^2}l d\eta\,.$$
Thus, if $f_l$ denotes a solution of \eqref{SteinEq}, then 
\begin{equation}\label{SteinEq1}
\E\left[\cN f_l(X_1,X_2)\right]=\int_{\R^2}\cN f_l(x_1,x_2)d\theta=\int_{\R^2}l d\theta-\int_{\R^2}ld\eta\,.
\end{equation}

The proof of the next lemma follows the well known proof for the one-dimensional case, and the only novelty is the upper bound  for the partial derivative with respect to the second variable.
\begin{lem}\label{SolStein0}
Let $l:\R^2\to\R$ be a continuously differentiable real valued function with bounded partial derivatives, $X_1\sim N(0,\sigma^2)$ and define 
$$f_l(x_1,x_2):=-\frac{1}{\sigma^2}\int_0^1\frac{1}{2\sqrt{t(1-t)}}\E\Big[X_1l\left(\sqrt{t}x_1+\sqrt{1-t}X_1,x_2\right)\Big]dt\,.$$
Then $f_l$ is a continuously differentiable real valued function such that: 
\begin{itemize}
\item $\|f_l\|_{\LL^\infty(\R^2)}\leq \|\partial_{x_1}l\|_{\LL^\infty(\R^2)}\,;$ 
\item $\|\partial_{x_1} f_l\|_{\LL^\infty(\R^2)}\leq\frac{1}{\sigma}\sqrt{\frac{2}{\pi}}\|\partial_{x_1} l\|_{\LL^\infty(\R^2)}\,;$
\item $\|\partial_{x_2} f_l\|_{\LL^\infty(\R^2)}\leq\frac{1}{\sigma}\sqrt{\frac{\pi}{2}}\|\partial_{x_2} l\|_{\LL^\infty(\R^2)}\,.$ 
\end{itemize}
Furthermore, $f_l$ is the unique bounded solution of \eqref{SteinEq}.
\end{lem}

\noindent\paragraph{\bf Proof of Lemma \ref{SolStein0}} 
Differentiating $f_l$ and carrying the derivative inside the integral and expectation can be justified using dominated convergence, and hence $f_l$ is a continuously differentiable function. In addition,    
\begin{equation}\label{repStein0}
\partial_{x_1} f_l(x_1,x_2)=-\frac{1}{\sigma^2}\int_0^1\frac{1}{2\sqrt{1-t}}\E\Big[X_1\partial_{x_1}l\left(\sqrt{t}x_1+\sqrt{1-t}X_1,x_2\right)\Big]dt\,.
\end{equation}
Fix $x_1$ and $x_2$ and consider the function $f(x)=l\left(\sqrt{t}x_1+\sqrt{1-t}x,x_2\right)$. Then
\begin{eqnarray*}
\E\Big[X_1l\left(\sqrt{t}x_1+\sqrt{1-t}X_1,x_2\right)\Big]&=&\E\Big[X_1f(X_1)\Big]\\
&=&\sigma^2\E\Big[f'(X_1)\Big]\\
&=&\sigma^2\sqrt{1-t} \E\Big[\partial_{x_1}l\left(\sqrt{t}x_1+\sqrt{1-t}X_1,x_2\right)\Big]\,,
\end{eqnarray*}
where we use Stein's lemma for the second equality, which leads to 
\begin{equation}\label{repStein1}
f_l(x_1,x_2)=-\int_0^1\frac{1}{2\sqrt{t}}\E\Big[\partial_{x_1}l\left(\sqrt{t}x_1+\sqrt{1-t}X_1,x_2\right)\Big]dt\,.
\end{equation}
Therefore,
\begin{eqnarray*}
\sigma^2 \partial_{x_1} f_l(x_1,x_2)- x_1 f_l(x_1,x_2)&=& \int_0^1\E\left[\left(-\frac{X_1}{2\sqrt{1-t}}+\frac{x_1}{2\sqrt{t}}\right)\partial_{x_1}l\left(\sqrt{t}x_1+\sqrt{1-t}X_1,x_2\right) \right]dt\\
&=&\int_0^1\E\left[\frac{d}{dt}l\left(\sqrt{t}x_1+\sqrt{1-t}X_1,x_2\right) \right]dt\\
&=&\E\left[\int_0^1\frac{d}{dt}l\left(\sqrt{t}x_1+\sqrt{1-t}X_1,x_2\right)dt \right]\\
&=&l(x_1,x_2)-\E\left[l(X_1,x_2)\right]\,,
\end{eqnarray*}
which shows that $f_l$ solves \eqref{SteinEq}. 
\medskip

Next we are going to use that 
$$\int_0^1\frac{1}{2\sqrt{t}}=\int_0^1\frac{1}{2\sqrt{1-t}} dt=1\mbox{ and }\int_0^1\frac{1}{2\sqrt{t(1-t)}}dt=\frac{\pi}{2}\,.$$
By \eqref{repStein1},
$$\| f_l\|_{\LL^\infty(\R^2)} \leq  \|\partial_{x_1}l\|_{\LL^\infty(\R^2)} \int_0^1\frac{1}{2\sqrt{t}} dt= \|\partial_{x_1}l\|_{\LL^\infty(\R^2)}\,,$$
and by \eqref{repStein0},
$$\|\partial_{x_1}f_l\|_{\LL^\infty(\R^2)} \leq  \frac{\E\left[|X_1|\right]}{\sigma^2}\|\partial_{x_1}l\|_{\LL^\infty(\R^2)}\int_0^1\frac{1}{2\sqrt{1-t}} dt=\frac{1}{\sigma}\sqrt{\frac{2}{\pi}} \|\partial_{x_1}l\|_{\LL^\infty(\R^2)}\,$$
(recall that $\E\left[|X_1|\right]=\sigma\sqrt{2\pi^{-1}}$). Now,
$$\partial_{x_2}f_l(x_1,x_2)=-\frac{1}{\sigma^2}\int_0^1\frac{1}{2\sqrt{t(1-t)}}\E\Big[X_1\partial_{x_2}l\left(\sqrt{t}x_1+\sqrt{1-t}X_1,x_2\right)\Big]dt\,,$$ 
and hence
\begin{eqnarray*}
\|\partial_{x_2}f_l\|_{\LL^\infty(\R^2)}&\leq &\frac{\E\left[|X_1|\right]}{\sigma^2}\|\partial_{x_2}l\|_{\LL^\infty(\R^2)}\int_0^1\frac{1}{2\sqrt{t(1-t)}}dt\\
&=&\frac{1}{\sigma}\sqrt{\frac{2}{\pi}}\|\partial_{x_2}l\|_{\LL^\infty(\R^2)} \frac{\pi}{2}\\
&=&\frac{1}{\sigma}\sqrt{\frac{\pi}{2}}\|\partial_{x_2}l\|_{\LL^\infty(\R^2)}\,.
\end{eqnarray*}
If $\tilde f$ is any other solution then 
$$\partial_{x_1}\left(e^{-\frac{{x_1}^2}{2\sigma^2}} \sigma^2\left(f_l(x_1,x_2)-\tilde f(x_1,x_2)\right)\right)=0 \,,$$
and thus $\tilde f(x_1,x_2)=f_l(x_1,x_2)+c(x_2)e^{\frac{x_1^2}{2\sigma^2}}$, which is bounded iff $c(x_2)\equiv 0$.

\hfill$\Box$\\ 

\begin{lem}\label{WassUpBound}
Let $\theta_t\equiv\P_{X_0^{\phi_1},X_t^{\phi_2}}$ and $\eta_t\equiv\P_{X_0^{\phi_1}}\otimes \P_{X_t^{\phi_2}}$. Then
$$\Wass(\eta_t,\theta_t)\leq \frac{1}{\|\phi_1\|_{\LL^2(\R)}}\sqrt{\frac{\pi}{2}} \E\Big[\left|\langle DX_t^{\phi_2},\phi_1\rangle_{\LL^2(\R)}\right|\Big]\,.$$ 
\end{lem}

\noindent\paragraph{\bf Proof Lemma \ref{WassUpBound}} 
By \eqref{IPCR} and \eqref{IP}, if $f$ is a continuously differentiable real valued function with bounded partial derivatives then $\E\left[W(\phi_1)f(X_0^{\phi_1},X_t^{\phi_2})\right]$ is given by  
$$\E\Big[\partial_{x_1} f(X_0^{\phi_1},X_t^{\phi_2})\langle DX_t^{\phi_1},\phi_1\rangle_{\LL^2(\R)}\Big]+\E\Big[\partial_{x_2} f(X_0^{\phi_1},X_t^{\phi_2})\langle DX_t^{\phi_2},\phi_1\rangle_{\LL^2(\R)}\Big]\,.$$
Recall that $X_0^{\phi_1}=\beta W(\phi_1)\sim N(0,\sigma^2)$, with $\sigma^2=\beta^2\big \|\phi_1\|^2_{\LL^2(\R)}$, and that $DX_0^{\phi_1}(z)=\beta\phi_1(z)$ (Lemma \ref{BasMalDer}). Thus, $\E\left[W(\phi_1)f(X_0^{\phi_1},X_t^{\phi_2})\right]$ is equal to 
$$\beta\big \|\phi_1\|^2_{\LL^2(\R)}\E\Big[\partial_{x_1} f(X_0^{\phi_1},X_t^{\phi_2})\Big]+\E\Big[\partial_{x_2} f(X_0^{\phi_1},X_t^{\phi_2})\langle DX_t^{\phi_2},\phi_1\rangle_{\LL^2(\R)}\Big]\,,$$
and hence
\begin{eqnarray*}
\E\left[X_0^{\phi_1} f(X^{\phi_1}_0,X_t^{\phi_2})\right]&=&\beta\E\left[W(\phi_1)f(X_0^{\phi_1},X_t^{\phi_2})\right]\\
&=&\sigma^2\E\Big[\partial_{x_1} f(X_0^{\phi_1},X_t^{\phi_2})\Big]+ \beta\E\Big[\partial_{x_2} f(X_0^{\phi_1},X_t^{\phi_2})\langle DX_t^{\phi_2},\phi_1\rangle_{\LL^2(\R)}\Big]\,,
\end{eqnarray*}
which implies that 
\begin{equation}\label{upbound}
\Big | \E\Big[\sigma^2\partial_{x_1} f(X_0^{\phi_1},X_t^{\phi_2})-X_0^{\phi_1} f(X_0^{\phi_1},X_t^{\phi_2})\Big]\Big |\leq \beta\|\partial_{x_2} f\|_{\LL^\infty(\R^2)} \E\Big[\left|\langle DX_t^{\phi_2},\phi_1\rangle_{\LL^2(\R)}\right |\Big]\,.
\end{equation}

By \eqref{SteinEq1}, Lemma \ref{SolStein0} and \eqref{upbound}, for every continuously differentiable real valued function with bounded partial derivatives $l:\R^2\to\R$, we have that
\begin{eqnarray}
\nonumber\Big|\int_{\R^2} ld\theta_t-\int_{\R^2} l d\eta_t\Big|&=&\Big|\E\left[\cN f_l\left(X_0^{\phi_1},X_t^{\phi_2}\right)\right]\Big|\\
\nonumber&=&\Big|\E\Big[\sigma^2\partial_{x_1} f_l(X_0^{\phi_1},X_t^{\phi_2})-X_0^{{\phi_1}}f_l(X_0^{\phi_1},X_t^{\phi_2})\Big]\Big|\\
\nonumber&\leq&\beta\|\partial_{x_2} f_l\|_{\LL^\infty(\R^2)}\E\Big[\left|\langle DX_t^{\phi_2},\phi_1\rangle_{\LL^2(\R)}\right|\Big]\\
\label{upbound1}&\leq &\frac{1}{\|\phi_1\|_{\LL^2(\R)}}\sqrt{\frac{\pi}{2}} \|\partial_{x_2} l\|_{\LL^\infty(\R^2)}\E\Big[\left|\langle DX_t^{\phi_2},\phi_1\rangle_{\LL^2(\R)}\right|\Big]\,.
\end{eqnarray}

The rest of the proof follows by approximating a Lipschitz function by continuously differentiable functions with bounded partial derivatives. Indeed, given $l\in\Lip_C$ and $\epsilon>0$ let  
$$l_\epsilon(x,y):=\E\left[l\left(x+\sqrt{\epsilon}N_1,y+\sqrt{\epsilon}N_2\right)\right]\,,$$
where $N_1$ and $N_2$ are independent standard normal random variables. Then $l_\epsilon$ is continuously differentiable with bounded partial derivatives. Furthermore,
$$\lim_{\epsilon\to 0}\|l_\epsilon-l\|_{\LL^\infty(\R^2)}=0\,\mbox{ and }\,\max\{ \|\partial_{x_1} l_{\epsilon}\|_{\LL^\infty(\R^2)}\,,\, \|\partial_{x_2} l_{\epsilon}\|_{\LL^\infty(\R^2)} \}\leq \|l_\epsilon\|_\Lip\leq \|l\|_\Lip\,.$$   
Recall that the Wasserstein distance \eqref{Wass} is defined by considering the collection of all Lipschitz functions $l$ such that $\|l\|_{\Lip}\leq 1$. Thus, if $\|l\|_{\Lip}\leq 1$ then $\|\partial_{x_2} l_\epsilon\|_{\LL^\infty(\R)}\leq 1$, and by \eqref{upbound1}
\begin{eqnarray*} 
\Big|\int_{\R^2} l d\theta_t-\int_{\R^2} l d\eta_t\Big| &\leq &\Big|\int_{\R^2} l_\epsilon d\theta_t-\int_{\R^2} l_\epsilon d\eta_t\Big| +2\|l_\epsilon-l\|_{\LL^\infty(\R^2)}\\
&\leq &\frac{1}{\|\phi_1\|_{\LL^2(\R)}}\sqrt{\frac{\pi}{2}} \E\Big[\left|\langle DX_t^{\phi_2},\phi_1\rangle_{\LL^2(\R)}\right|\Big]+2\|l_\epsilon-l\|_{\LL^\infty(\R^2)}\,.
\end{eqnarray*}
By letting $\epsilon\to 0$, we get that 
$$\Big|\int_{\R^2} l d\theta_t-\int_{\R^2} l d\eta_t\Big|\leq \frac{1}{\|\phi_1\|_{\LL^2(\R)}}\sqrt{\frac{\pi}{2}}\E\Big[\left|\langle DX_t^{\phi_2},\phi_1\rangle_{\LL^2(\R)}\right|\Big]\,,$$
which concludes the proof of the lemma (using Remark \ref{Correction}).

\hfill$\Box$\\

\noindent\paragraph{\bf Proof of Theorem \ref{AsyInd}} It is a direct consequence of Lemma \ref{ConvFor}, Lemma \ref{ExiMalDer} and Lemma \ref{WassUpBound}.

\hfill$\Box$\\

\end{document}